\documentclass[smallextended,numbook,runningheads]{svjour3}		
\smartqed	
\usepackage{graphicx}
\usepackage{amsmath}
\usepackage{amsfonts,graphics,amssymb,sidecap}
\usepackage{ntheorem}
\usepackage{epstopdf}	
\usepackage[mathscr]{euscript}
\usepackage{times}
\usepackage{color}
\DeclareMathOperator{\diam}{diam}			
\DeclareMathOperator{\dist}{dist}			
\DeclareMathOperator{\lcm}{lcm}				
\DeclareMathOperator{\lip}{lip}				
\newcommand{\C}{{\mathbb C}}
\newcommand{\N}{{\mathbb N}}
\newcommand{\R}{{\mathbb R}}
\renewcommand{\S}{{\mathbb S}}
\newcommand{\Z}{{\mathbb Z}}
\newcommand{\cA}{{\mathcal A}}
\newcommand{\cB}{{\mathcal B}}
\newcommand{\cV}{{\mathcal V}}
\newcommand{\cW}{{\mathcal W}}
\newcommand{\sF}{{\mathscr F}}

\newcommand{\sH}{{\mathscr H}}
\newcommand{\sK}{{\mathscr K}}
\newcommand{\sX}{{\mathscr X}}
\newcommand{\sY}{{\mathscr Y}}
\renewcommand{\d}{\,{\mathrm d}}
\newcommand{\eps}{\varepsilon}
\newcommand{\tm}{\times}
\newcommand{\vphi}{\varphi}

\newcommand{\intoo}[1]{\left(#1\right)}				
\newcommand{\intcc}[1]{\left[#1\right]}				
\newcommand{\set}[1]{\left\{#1\right\}}				
\newcommand{\abs}[1]{\left|#1\right|}				
\newcommand{\norm}[1]{\left\|#1\right\|}			
\newcommand{\fall}{\quad\text{for all }}
\newcommand{\cref}[1]{Cor.~\ref{#1}}

\newcommand{\fref}[1]{Fig.~\ref{#1}}
\newcommand{\href}[1]{Hyp.~\ref{#1}}
\newcommand{\lref}[1]{Lemma~\ref{#1}}
\newcommand{\pref}[1]{Prop.~\ref{#1}}

\newcommand{\sref}[1]{Sec.~\ref{#1}}
\newcommand{\tref}[1]{Thm.~\ref{#1}}
\newcommand{\Cd}{C_d^0}
\newcommand{\Cda}{C_d^\alpha}

\numberwithin{equation}{section}
\newtheorem{hypothesis}{Hypothesis}[section]
\allowdisplaybreaks
\allowdisplaybreaks
\begin{document}
\title{Numerical dynamics of integrodifference equations}
\subtitle{Periodic solutions and invariant manifolds in $C^\alpha(\Omega)$}
\titlerunning{Numerical dynamics of integrodifference equations}	
\author{Christian P\"otzsche}
%
%
\institute{Christian P\"otzsche\at 
Institut f\"ur Mathematik, 
Universit\"at Klagenfurt, 
Universit\"atsstra{\ss}e 65--67, 
A-9020 Klagenfurt, Austria, 
\email{christian.poetzsche@aau.at}}

\date{Received: \today / Accepted: date}
%
\maketitle
\begin{abstract}
	Integrodifference equations are versatile models in theoretical ecology for the spatial dispersal of species evolving in non-overlapping generations. The dynamics of these infinite-dimensional discrete dynamical systems is often illustrated using computational simulations. This paper studies the effect of Nystr\"om discretization to the local dynamics of periodic integrodifference equations with H\"older continuous functions over a compact domain as state space. We prove persistence and convergence for hyperbolic periodic solutions and their associated stable and unstable manifolds respecting the convergence order of the quadrature/cubature method. 


\end{abstract}
\section{Introduction}
Integrodifference equations (IDEs for short) are infinite-dimensional discrete-time dynamical systems. They became popular tools in theoretical ecology over the recent years modelling the temporal evolution and spatial dispersal of species having non-overlapping generations \cite{kot:schaeffer:86,lutscher:19}. Furthermore, IDEs canonically arise as time discretizations of integrodifferential equations, as time-$1$-map of evolutionary partial differential equations or in the iterative solution of (nonlinear) boundary value problems \cite[p.~190]{martin:76}. It is understood that IDEs involve an integral operator which is typically of Hammerstein- or more general Urysohn-type. Indeed, for our purposes a sufficiently flexible class are recursions of the form
\begin{equation}
	\tag{$I_0$}
	u_{t+1}(x)=\int_\Omega f_t(x,y,u_t(y))\d y\fall x\in\Omega, 
	\label{ide0}
\end{equation}
whose natural state spaces consists of continuous or integrable functions over a compact subset $\Omega\subset\R^\kappa$ (the habitat in applications from ecology). In applied sciences the long-term behavior of IDEs is willingly illustrated using numerical simulations. For this purpose, \cite[pp.~112--113]{lutscher:19} suggests to replace the integral in \eqref{ide0} by the trapezoidal or the Simpson rule. Both are special cases of general Nystr\"om methods
\begin{equation}
	\tag{$I_n$}
	u_{t+1}(x)=\sum_{\eta\in\Omega_n}w_\eta f_t(x,\eta,u_t(\eta))\fall x\in\Omega
	\label{iden}
\end{equation}
based on convergent quadrature or cubature rules with weights $w_\eta\geq 0$ and nodes $\eta\in\Omega_n$ over a finite grid $\Omega_n\subset\Omega$. Here, $n\in\N$ is related to the number of nodes in $\Omega_n$ and therefore the accuracy of the approximation, see \cite{davis:rabinowitz:84}. We point out that Nystr\"om methods yield full discretizations of \eqref{ide0} and can be evaluated immediately. 

While the numerical analysis of integral equations is a well-established field, e.g.\ \cite{atkinson:92,hackbusch:95}, this paper enriches it by a dynamical aspect: We study and relate the long-term behavior of the iterates $u_t$ generated by an IDE \eqref{ide0} to those of a Nystr\"om discretization \eqref{iden}. This brings us to the area of numerical dynamics \cite{lubich:01,stuart:95,stuart:humphries:98} addressing the following questions:
\begin{itemize}
	\item Which dynamical or asymptotic properties of an IDE \eqref{ide0} as $t\to\infty$ are preserved when passing to its Nystr\"om discretizations \eqref{iden}? 

	\item What can be said about convergence as $n\to\infty$ when the approximations become increasingly more accurate? In particular, are convergence rates of the integration rules preserved? 
\end{itemize}
Concerning the qualitative behavior of autonomous ODEs, such problems originate in \cite{beyn:87} and are surveyed in \cite{stuart:humphries:98}. In between various contributions to continuous-time infinite-di\-men\-sional dynamical systems generated by functional differential equations \cite{garay:05} or evolutionary (e.g.\ parabolic) partial differential equations \cite{stuart:95,lubich:01} arose, both for spatial, as well as for full discretizations. IDEs merely require spatial discretization, but have in common with these problems that conventional error estimates fail to describe asymptotic behavior. In fact, bounds for the global discretization error typically grow exponentially in time and therefore establish convergence as $n\to\infty$ only over compact time intervals \cite[Thm.~4.1]{poetzsche:18a}. Thus, techniques extending those of classical numerical analysis are required to tackle the above problems. 

Previous contributions to the numerical dynamics of IDEs address basics and error estimates \cite{poetzsche:18a}, as well as the persistence/convergence of globally asymptotically stable solutions \cite{poetzsche:18b}. This paper focusses on an another important aspect, namely the local saddle-point structure near periodic solutions to \eqref{ide0}. Related work, but for autonomous evolutionary differential equations near equilibria, is due to \cite{beyn:87,garay:93} (ODEs), \cite{alouges:debussche:91} (parabolic PDEs) and \cite{farkas:01} (retarded FDEs). 

In contrast, we study time-periodic IDEs \eqref{ide0} in the vicinity of periodic solutions. We stress that periodic time-dependence is strongly motivated by applications to incorporate seasonal influences. While \cite{poetzsche:18a,poetzsche:18b} apply to semi-dis\-cre\-ti\-zations of \eqref{ide0} of collocation- or degenerate kernel-type \cite{atkinson:92,hackbusch:95}, which act between finite-dimensional function spaces, but still contain integrals, we tackle Nystr\"om discretizations \eqref{iden}, because they can be evaluated immediately. At this point the question for an ambient state space of \eqref{ide0} arises. A natural choice are the continuous functions $C(\Omega)$ over a compact $\Omega\subset\R^\kappa$. Here however, already for linear integral operators, the discretization error under Nystr\"om methods converges only in the strong, but not in the uniform topology as $n\to\infty$, see \cite[pp.~130--131, Lemma~4.7.6]{hackbusch:95}. Using the theory of collectively compact operators \cite{anselone:71} one can still establish that fixed-points to \eqref{ide0} (and their stability properties) persist \cite{atkinson:73,weiss:74}. Nonetheless, it is not clear how to establish convergence of the associated stable and unstable manifolds of \eqref{iden} to those of the original problem \eqref{ide0}. For this reason we retreat to the H\"older continuous functions $C^\alpha(\Omega)$ as state space. This set-up is sufficiently general to capture most relevant applied problems \cite{kot:schaeffer:86,lutscher:19} and has the advantage that a more conventional perturbation theory (see App.~\ref{appA}) applies to realize our goals. It should not be concealed, though, that the prize for this endeavor are more involved assumptions and technical preliminaries on Urysohn operators (well-definedness, complete continuity, differentiability). For the sake of a brief presentation they are outsourced to \cite{poetzsche:20,poetzsche:20b}. 

The structure of our presentation is as follows: In \sref{sec2} we introduce a flexible framework for general periodic difference equations in Banach spaces and their linearization. Perturbation results for the Floquet spectrum of linear periodic equations are given in \tref{thmhyper}, while \tref{thmsol} addresses persistence and convergence of hyperbolic solutions and \tref{thmbundle} the associated stable and unstable manifolds -- when dealing with periodic equations we speak of fiber bundles. Although tailor-made for Nystr\"om discretizations of IDEs, these results also apply to collocation- or degenerate kernel-discretizations, as well as when studying time-periodic evolutionary differential equations via their time-$h$-maps. The concrete case of Urysohn IDEs \eqref{ide0} is saved for \sref{sec3} and illustrates how Thms.~\ref{thmhyper}--\ref{thmbundle} apply. One obtains convergence of both hyperbolic solutions, and of the functions parametrizing their invariant fiber bundles with a rate given by the H\"older exponent $\alpha\in(0,1]$ of the kernel functions $f_t$ in the first variable. Nevertheless, for smooth $f_t$ the higher-order convergence rates inherited from the particular quadrature/cubature rules are established. 
\paragraph{Notation}
We write $\R_+:=[0,\infty)$ for the nonnegative reals, $\S^1:=\set{z\in\C:\,\abs{z}=1}$ for the unit circle in $\C$, $[\cdot]:\R\to\Z$ is the integer function and $\abs{\cdot}$ denotes norms on finite-dimen\-sional spaces. On the Cartesian product $X\tm Y$ of normed spaces $X,Y$, 
\begin{equation}
	\norm{(x,y)}:=\max\bigl\{\norm{x}_X,\norm{y}_Y\bigr\}
	\label{prodnorm}
\end{equation}
is the product norm and we proceed accordingly on products of more than two spaces. The open resp.\ closed balls in $X$ with radius $r\geq 0$ and center $x\in X$ are
\begin{align*}
	B_r(x,X)&:=\set{y\in X:\,\norm{y-x}<r},&
	\bar B_r(x,X)&:=\set{y\in X:\,\norm{y-x}\leq r};
\end{align*}
on a finite-dimensional $X$ we write $B_r(x)$ and $\bar B_r(x)$. For nonempty $A\subseteq X$, $\diam A$ denotes the diameter of $A$, $\dist_A(x):=\inf_{a\in A}\norm{x-a}$ the distance of a point $x\in X$ from $A$, $\dist(B,A):=\sup_{b\in B}\dist_A(b)$ the Hausdorff semidistance of $B\subseteq X$ from $A$ and we set $B_r(A):=\{x\in X:\,\dist_A(x)<r\}$. We denote a subset $\cA\subseteq\Z\tm X$ as \emph{nonautonomous set} with fibers $\cA(t):=\{x\in X:\,(t,x)\in\cA\}$, $t\in\Z$ and write 
$
	\cB_r(\phi):=\set{(t,u)\in\Z\tm X:\,\norm{u-\phi_t}<r}
$
for the $r$-neighborhood of a sequence $\phi=(\phi_t)_{t\in\Z}$ in $X$. 

The bounded $k$-linear maps from the Cartesian product $X^k$ to $Y$ are denoted by $L_k(X,Y)$, $L(X,Y):=L_1(X,Y)$ and $L_0(X,Y):=Y$. Moreover, we abbreviate $L_k(X):=L_k(X,X)$, $L(X):=L(X,X)$, $GL(X)$ are the invertible maps in $L(X)$ and $I_X$ is the identity on $X$. Furthermore, $N(T)$ is the kernel and $R(T)$ the range of $T\in L(X,Y)$; $\sigma(S)$ is the spectrum and $\sigma_p(S)$ the point spectrum of $S\in L(X)$. 
\section{Difference equations and perturbation}
\label{sec2}
Let $(X,\norm{\cdot})$ denote a Banach space. 
\subsection{Periodic difference equations}
Abstractly, we are interested in a family of nonautonomous difference equations
\begin{equation}
	\tag{$\Delta_n$}
	\boxed{u_{t+1}=\sF_t^n(u_t)}
	\label{deq}
\end{equation}
with right-hand sides $\sF_t^n:U_t\to X$ on open sets $U_t\subseteq X$, $t\in\Z$, parametrized by $n\in\N_0$. In the following, $n\in\N$ is a discretization parameter such that $\sF_t^n$ are understood as approximations converging to the original problem $\sF_t^0$ as $n\to\infty$ in a sense to be defined below. A nonautonomous set $\cA\subseteq\Z\tm X$ with fibers $\cA(t)\subseteq U_t$ for all $t\in\Z$ is called \emph{forward invariant} or \emph{invariant} (w.r.t.\ \eqref{deq}), provided \begin{align*}
	\sF_t^n(\cA(t))&\subseteq\cA(t+1),&
	\sF_t^n(\cA(t))&=\cA(t+1)\fall t\in\Z
\end{align*}
resp., holds. Given an initial time $\tau\in\Z$, a \emph{forward solution} to \eqref{deq} is a sequence $\phi=(\phi_t)_{\tau\leq t}$ satisfying $\phi_t\in U_t$ and the \emph{solution identity} $\phi_{t+1}=\sF_t^n(\phi_t)$ for all $\tau\leq t$, a \emph{backward solution} fulfills the solution identity for $t<\tau$ and for an \emph{entire solution} $(\phi_t)_{t\in\Z}$ one has $\phi_{t+1}\equiv\sF_t^n(\phi_t)$ on $\Z$. The forward solution starting at $\tau$ in the initial state $u_\tau\in U_\tau$ is uniquely determined as composition
$$
	\vphi^n(t;\tau,u_\tau)
	:=
	\begin{cases}
		\sF_{t-1}^n\circ\ldots\circ\sF_\tau^n(u_\tau),&\tau<t,\\
		u_\tau,&t=\tau
	\end{cases}
$$
and denoted as the \emph{general solution} to \eqref{deq}; it is defined as long as the compositions stay in $U_t$. A difference equation \eqref{deq} is called $\theta_0$-\emph{periodic}, if both $\sF_{t+\theta_0}^n=\sF_t^n$ and $U_{t+\theta_0}=U_t$ hold for all $t\in\Z$ with some \emph{basic period} $\theta_0\in\N$; an \emph{autonomous} equation is $1$-periodic, i.e.\ there exists a $\sF^n:U\to X$ with $\sF_t^n\equiv\sF^n$, $U_t\equiv U$ on $\Z$. A $\theta_1$-\emph{periodic solution} to \eqref{deq} is an entire solution satisfying $\phi_t\equiv\phi_{t+\theta_1}$ on $\Z$. 

Given a fixed $\theta\in\N$ and a sequence $u=(u_t)_{t\in\Z}$ with $u_t\in U_t$, $t\in\Z$, let us introduce the open product $\hat U:=U_0\tm\ldots\tm U_{\theta-1}$ and
\begin{align*}
	\hat u&:=(u_0,\ldots,u_{\theta-1})\in\hat U,&
	(\overline{u_0,\ldots,u_{\theta-1}})&:=(u_{t\text{ mod }\theta})_{t\in\Z}.
\end{align*}
In order to characterize and compute periodic solutions to \eqref{deq}, $n\in\N_0$, we introduce the nonlinear operators
\begin{align}
	\hat\sF^n:\hat U&\to X^\theta,&
	\hat\sF^n(\hat u)
	&:=
	\begin{pmatrix}
		\sF_{\theta-1}^n(u_{\theta-1})\\
		\sF_0^n(u_0)\\
		\sF_1^n(u_1)\\
		\vdots\\
		\sF_{\theta-2}^n(u_{\theta-2})
	\end{pmatrix}
	\label{fhatdef}
\end{align}
and use the norm induced via \eqref{prodnorm} on the Cartesian product $X^\theta$. 

The next two results are immediate: 
\begin{lemma}\label{lemlinb}
	Let $n\in\N_0$, \eqref{deq} be $\theta_0$-periodic and $\theta$ be a multiple of $\theta_0$: 
	\begin{enumerate}
		\item If $(\phi_t)_{t\in\Z}$ is a $\theta$-periodic solution to \eqref{deq}, then $\hat\phi\in\hat U$ is a fixed point of $\hat\sF^n$. 

		\item Conversely, if $\hat\phi\in\hat U$ is a fixed point of $\hat\sF^n$, then $(\overline{\phi_0,\ldots,\phi_{\theta-1}})$ is a $\theta$-periodic solution to \eqref{deq}. 
	\end{enumerate}
\end{lemma}
This characterization of periodic solutions to \eqref{deq} via the mapping $\hat\sF$ has the numerical advantage to avoid the computation of compositions $\vphi^n(\theta+\tau;\tau,\cdot):U_\tau\to X$, $\tau\in\Z$, and therefore preserves (numerical) backward stability (see \cite{fassbender:kresser:06}). 
\begin{lemma}\label{lemhatder}
	Let $n\in\N_0$, $m\in\N$, \eqref{deq} be $\theta_0$-periodic and $\theta$ be a multiple of $\theta_0$. If every $\sF_t^n:U_t\to X$, $0\leq t<\theta_0$, is $m$-times continuously (Fr{\'e}chet) differentiable, then $\hat\sF^n:\hat U\to X^\theta$ is of class $C^m$ and for every $\hat u\in\hat U$ one has
	$$
		D\hat\sF^n(\hat u)
		=
		\begin{pmatrix}
			0 & 0 & \cdots & \cdots & D\sF_{\theta-1}^n(u_{\theta-1})\\
			D\sF_0^n(u_0) & 0 & \cdots & \cdots & 0\\
			0 & D\sF_1^n(u_1) & 0 & \cdots & 0\\
			\vdots & \ddots & \ddots & \ddots & \vdots\\
			0 & \cdots & 0 & D\sF_{\theta-2}^n(u_{\theta-2}) & 0
		\end{pmatrix}. 
	$$
\end{lemma}
\subsection{Linear periodic difference equations}
\label{sec22}
Suppose that $\sK_t^n\in L(X)$, $t\in\Z$, and consider a family of linear difference equations
\begin{equation}
	\tag{$L_n$}
	\boxed{u_{t+1}=\sK_t^nu_t}
	\label{lin}
\end{equation}
in $X$ parametrized by $n\in\N_0$. As above we understand \eqref{lin}, $n\in\N$, as perturbations of an initial problem $(L_0)$. The \emph{transition operator} $\Phi^n(t,\tau)\in L(X)$ of \eqref{lin} is 
\begin{equation}
	\Phi^n(t,\tau):=
	\begin{cases}
		\sK_{t-1}^n\cdots\sK_\tau^n,&\tau<t,\\
		I_X,&t=\tau.
	\end{cases}
	\label{trans}
\end{equation}
We are interested in $\theta$-periodic equations \eqref{lin}, that is
\begin{equation}
	\sK_t^n=\sK_{t+\theta}^n\fall t\in\Z,
	\label{kper}
\end{equation}
allowing us to introduce the \emph{period operator} $\Xi_\theta^n:=\Phi^n(\theta,0)\in L(X)$ of \eqref{lin}. Its eigenvalues are the \emph{Floquet multipliers} and $\sigma_p(\Xi_\theta^n)$ is the \emph{Floquet spectrum} of \eqref{lin}. 

One says a linear difference equation \eqref{lin} is \emph{weakly hyperbolic}, if $1\not\in\sigma(\Xi_\theta^n)$, and \emph{hyperbolic}, if $\S^1\cap\sigma(\Xi_\theta^n)=\emptyset$ holds. In the hyperbolic situation, the spectrum can be decomposed as $\sigma(\Xi_\tau^n)=\sigma_+\dot\cup\sigma_-$ with spectral sets
\begin{align*}
	\sigma_+&\subseteq B_1(0),&
	\sigma_-&\subseteq\C\setminus\bar B_1(0). 
\end{align*}
With the spectral projections 
$P_+^n:=\frac{1}{2\pi\i}\int_{\S^1}(zI_{X}-\Xi_\theta^n)^{-1}\d z$, 
$P_-^n:=I_X-P_+^n$ 
we introduce the fibers $\cV_+^n(t):=\Phi^n(t,0)R(P_+^n)$ and $\cV_-^n(t):=\Phi^n(t,0)R(P_-^n)$, 
first for $t\geq 0$ and then by $\theta$-periodic continuation on $\Z$. This yields $\theta$-periodic nonautonomous sets $\cV_+^n\subseteq\Z\tm X$ (\emph{stable vector bundle}) and $\cV_-^n\subseteq\Z\tm X$ (\emph{unstable vector bundle}) of \eqref{lin}. Then $\cV_+^n$ is forward invariant, while $\cV_-^n$ is invariant w.r.t.\ \eqref{lin}.

For compact operators $\Xi_\theta^n\in L(X)$ the Riesz-Schauder theory \cite[pp.~428ff]{lang:93} applies: Every Floquet multiplier $\lambda\in\sigma_p(\Xi_\theta^n)$ possesses a minimal $\iota(\lambda)\in\N$ so that $N(\lambda I_X-\Xi_\theta^n)^j=N(\lambda I_X-\Xi_\theta^n)^{j+1}$ for all $j\geq\iota(\lambda)$ leading to finite-dimensional \emph{generalized eigenspaces} $N(\lambda I_X-\Xi_\theta^n)^{\iota(\lambda)}$. All unstable fibers $\cV_-^n(t)$, $t\in\Z$, have a constant finite dimension, which is denoted as the \emph{Morse index} of \eqref{lin} and equals the finite sum $\sum_{\lambda\in\sigma_-}\dim N(\lambda I_X-\Xi_\theta^n)^{\iota(\lambda)}$ of algebraic multiplicities. 

We begin with a perturbation result for hyperbolic linear systems \eqref{lin} under uniform convergence: 
\begin{theorem}[perturbed hyperbolicity]\label{thmhyper}
	Suppose that the $\theta$-periodic linear difference equations \eqref{lin}, $n\in\N_0$, fulfill: 
	\begin{itemize}
		\item[(i)] $\lim_{n\to\infty}\norm{\sK_t^n-\sK_t^0}_{L(X)}=0$ for all $0\leq t<\theta$, 

		\item[(ii)] $\Xi_\theta^n\in L(X)$ is compact for all $n\in\N$. 
	\end{itemize}
	Then also the period operator $\Xi_\theta^0\in L(X)$ of $(L_0)$ is compact and there exists a $N\in\N$ such that the following holds for all $n\geq N$ or $n=0$:
	\begin{enumerate}
		\item With $(L_0)$ also the perturbed equation \eqref{lin} is weakly hyperbolic,

		\item with $(L_0)$ also the perturbed equation \eqref{lin} is hyperbolic. In particular, for reals $\beta\in\bigl(\max\set{0,1-\tfrac{1}{2}\dist(\sigma(\Xi_\theta^0),\S^1)},1\bigr)$, there exists a $\theta$-periodic sequence $(P_t^n)_{t\in\Z}$ of \emph{invariant projectors} in $L(X)$ with $\sK_t^nP_t^n=P_{t+1}^n\sK_t^n$ for all $t\in\Z$, so that the transition operators $\Phi^n(t,s)$ satisfy $\dim\cV_-^n=\dim\cV_-^0$ and the estimates
		\begin{equation}
			\begin{split}
				\norm{\Phi^n(t,s)P_s^n}_{L(X)}&\leq K\beta^{t-s}\fall s\leq t,\\
				\norm{\Phi^n(t,s)[I_X-P_s^n]}_{L(X)}&\leq K\beta^{s-t}\fall t\leq s,
			\end{split}
			\label{thmhyper2}
		\end{equation}

		\item $\lim_{n\to\infty}\norm{P_t^n-P_t^0}_{L(X)}=0$ for all $t\in\Z$.
	\end{enumerate}
\end{theorem}
\begin{proof}
	Let $0\leq s<\theta$. Due to (i) the sequence $\bigl(\norm{\sK_s^n-\sK_s^0}\bigr)_{n\in\N}$ is bounded and consequently we obtain from $\norm{\sK_s^n}\leq\norm{\sK_s^0}+\norm{\sK_s^n-\sK_s^0}$ and the periodicity condition \eqref{kper} that $c_t:=\sup_{n\in\N_0}\norm{\sK_t^n}<\infty$ for all $t\in\Z$. 

	(I) \emph{Claim}: $\lim_{n\to\infty}\|\Phi^n(t,0)-\Phi^0(t,0)\|=0$ for all $0\leq t$.\\
	We proceed by mathematical induction. Thanks to \eqref{trans}, for $t=0$ the assertion is trivial and for $t=1$ it results from (i). In the induction step $t\to t+1$ we obtain
	\begin{align*}
		&
		\|\Phi^n(t+1,0)-\Phi^0(t+1,0)\|
		\stackrel{\eqref{trans}}{=}
		\|\sK_t^n\Phi^n(t,0)-\sK_t^0\Phi^0(t,0)\|\\
		&\leq
		\|\sK_t^n\|\|\Phi^n(t,0)-\Phi^0(t,0)\|+
		\|\sK_t^n-\sK_t^0\|\|\Phi^0(t,0)\|\\
		&\leq
		c_t\|\Phi^n(t,0)-\Phi^0(t,0)\|
		+
		\|\sK_t^n-\sK_t^0\|\prod_{r=0}^{t-1}c_r
		\xrightarrow[n\to\infty]{(i)}0
	\end{align*}
	from the induction hypothesis and the triangle inequality, yielding the claim. 

	(II) \emph{Claim}: $\Xi_\theta^0\in L(X)$ is compact.\\
	If we set $t=\theta$ in claim (I), then the period operators satisfy 
	\begin{equation}
		\lim_{n\to\infty}\norm{\Xi_\theta^n-\Xi_\theta^0}=0.
		\label{limit}
	\end{equation}
	Hence, $\Xi_\theta^0$ is the uniform limit of by (ii) compact operators $\Xi_\theta^n$, $n\in\N$, and consequently compact \cite[p.~416, Thm.~1.1]{lang:93}. 

	(III) \emph{Claim}: For every nonempty closed $S\subseteq\C$ with $\sigma(\Xi_\theta^0)\cap S=\emptyset$ there exists a $n_1\in\N$ such that $\sigma(\Xi_\theta^n)\cap S=\emptyset$ for all $n\geq n_1$.\\
	Since the closed $S$ and the compact $\sigma(\Xi_\theta^0)$ are disjoint they have a positive distance. Therefore, there is an $\eps>0$ so that $S\cap B_\eps(\sigma(\Xi_\theta^0))=\emptyset$. By the upper semicontinuity of the spectrum \cite[p.~80, Lemma~3]{baumgartel:85} and relation \eqref{limit} there is a $n_1\in\N$ with $\sigma(\Xi_\theta^n)\subset B_\eps(\sigma(\Xi_\theta^0))$ and consequently $\sigma(\Xi_\theta^n)$ stays disjoint from $S$ for all $n\geq n_1$. 

	(a) If $(L_0)$ is weakly hyperbolic, then $\sigma(\Xi_\theta^0)\cap\set{1}=\emptyset$ and (III) applied to the singleton $S=\set{1}$ yields the assertion. 

	(b) The hyperbolicity of \eqref{lin} results as above in (a) with $S=\S^1$. Furthermore, then \cite[p.~44, Prop.~3.13]{russ:15} implies that \eqref{lin} possess an exponential dichotomy on $\Z$ as claimed with the $\theta$-periodic invariant projectors $P_t^n$ satisfying 
	\begin{equation}
		\sK_t^n|_{N(P_t^n)}\in GL(N(P_t^n),N(P_{t+1}^n))\fall t\in\Z
		\label{inv}
	\end{equation}
	and $I_X-P_0^n=P_-^n$. In particular, by \eqref{limit} and \cite[p.~80, Cor.~1]{baumgartel:85} the spectral projections $P_-^n$ associated to the unstable spectral parts of $\Xi_\theta^n$, $n\in\N_0$, fulfill that $\dim R(P_-^n)=\dim R(P_-^0)$ for large $n$, say for $n\geq n_2$. Thanks to \eqref{inv} this extends to the dimension of the unstable bundles $\cV_-^n$. Finally, we set $N:=\max\set{n_1,n_2}$. 

	(c) Combining \eqref{limit} with \cite[p.~80, Lemma~4]{baumgartel:85} yields that the spectral projections satisfy $\lim_{n\to\infty}\norm{P_-^n-P_-^0}=0$. Together with claim (I) we obtain for $t\in\Z$ that
	\begin{align*}
		&
		\norm{P_t^n-P_t^0}
		=
		\norm{[I_X-P_t^n]-[I_X-P_t^0]}\\
		\stackrel{\eqref{inv}}{\leq}&
		\norm{\Phi^n(t,0)}\norm{P_-^n-P_-^0}\norm{\Phi^0(0,t)}
		+
		\norm{\Phi^n(t,0)-\Phi^0(t,0)}\norm{P_-^0\Phi^0(0,t)}\\
		\leq&
		\bigg(\prod_{r=0}^{t-1}c_r\bigg)\norm{P_-^n-P_-^0}\norm{\Phi^0(0,t)}
		+
		\norm{\Phi^n(t,0)-\Phi^0(t,0)}\norm{P_-^0\Phi^0(0,t)}
	\end{align*}
	from the triangle inequality, whose right-hand side converges to $0$ as $n\to\infty$. 
	\qed
\end{proof}
\subsection{Perturbation of hyperbolic solutions and invariant bundles}
We next address the robustness of $\theta_1$-periodic solutions $\phi^0$ to general $\theta_0$-periodic difference equations $(\Delta_0)$, as well as their nearby saddle-point structure consisting of stable and unstable bundles (see \cite[pp.~143ff, Chap.~6]{irwin:80}, \cite[pp.~256ff, Sect.~4.6]{poetzsche:10b}) under perturbation. By imposing a natural hyperbolicity condition on the solution $\phi^0$ it is shown that also the perturbations \eqref{deq} have (locally unique) periodic solutions $\phi^n$ for sufficiently large $n$, which converge to $\phi^0$ in the limit $n\to\infty$. 

Let $\theta:=\lcm\set{\theta_0,\theta_1}$. We suppose that the right-hand sides $\sF_t^n$ of \eqref{deq} are continuously differentiable. Our endeavor is based on the \emph{variational equations}
\begin{equation}
	\tag{$V_n$}
	\boxed{v_{t+1}=D\sF_t^n(\phi_t^n)v_t}
	\label{var}
\end{equation}
associated to $\theta$-periodic solutions $\phi^n$ of \eqref{deq}, $n\in\N_0$. Since the linear equations \eqref{var} are $\theta$-periodic, the terminology and results from \sref{sec22} apply to \eqref{var} with $\sK_t^n=D\sF_t^n(\phi_t^n)$ and the period operator $\Xi_\theta^n$, $n\in\N_0$. In this context, we understand a solution $\phi^n$ of \eqref{deq} as \emph{(weakly) hyperbolic}, if \eqref{var} has the corresponding property. 
\begin{lemma}\label{lemhatdera}
	Let $n\in\N_0$. If $\sF_t^n:U_t\to X$ is continuously differentiable for all $0\leq t<\theta_0$, then the derivatives of the mappings $\hat\sF^n:\hat U\to X^\theta$ defined in \eqref{fhatdef} satisfy $\sigma(\Xi_\theta^n)\setminus\set{0}=\sigma(D\hat\sF^n(\hat\phi^n))^\theta\setminus\set{0}$ and $\sigma_p(\Xi_\theta^n)\setminus\set{0}=\sigma_p(D\hat\sF^n(\hat\phi^n))^\theta\setminus\set{0}$.
\end{lemma}
Based on this result, the (Floquet) spectrum of \eqref{var} can be computed from the (point) spectrum of the cyclic block operator given in \lref{lemhatder}. This has the numerical advantage of avoiding to evaluate the compositions (matrix products) $\Xi_\theta^n$. 
\begin{proof}
	Keeping $n\in\N_0$ fixed, we abbreviate $\sK_t=D\sF_t^n(\phi_t^n)$, $t\in\Z$, and observe that the $\theta$th power of $D\hat\sF(\hat\phi^n)$ given in \lref{lemhatder} becomes a block diagonal operator
	$$
		D\hat\sF(\hat\phi^n)^\theta
		=
		\begin{pmatrix}
			\sK_{\theta-1}\sK_{\theta-2}\cdots\sK_0 & & &\\
			& \sK_0\sK_{\theta-1}\cdots\sK_1 & &\\
			& & \ddots &\\
			& & & \sK_{\theta-2}\cdots\sK_0\sK_{\theta-1}
		\end{pmatrix}.
	$$
	Referring to \cite[p.~42, Prop.~3.11(a)]{russ:15} one has 
	$
		\sigma(\Xi_\theta^n)\setminus\set{0}
		=
		\sigma(\sK_{t+\theta-1}\cdots\sK_t)\setminus\set{0}
	$
	for all $t\in\Z$ and therefore $\sigma(\Xi_\theta^n)\setminus\set{0}=\sigma(D\hat\sF(\hat\phi^n)^\theta)\setminus\set{0}$. Now the Spectral Mapping Theorem \cite[p.~65, Thm.~2]{baumgartel:85} yields the assertion for the spectra. Concerning the point spectrum the claim follows directly from the corresponding eigenvalue-eigenvector relations and the solution identity for \eqref{var}. 
	\qed
\end{proof}

Our next result establishes persistence of hyperbolic periodic solutions to $(\Delta_0)$:
\begin{theorem}[perturbed periodic solutions]\label{thmsol}
	Let $\theta=\lcm\set{\theta_0,\theta_1}$. Suppose that the $\theta_0$-periodic difference equations \eqref{deq}, $n\in\N_0$, fulfill: 
	\begin{itemize}
		\item[(i)] $\sF_t^n:U_t\to X$ are continuously differentiable for all $0\leq t<\theta_0$ and $n\in\N_0$, 

		\item[(ii)] $D\sF_t^n:U_t\to L(X)$, $n\in\N$, are uniformly continuous on bounded sets uniformly in $n\in\N$, the family $\set{D\sF_t^n}_{n\in\N}$ is equicontinuous for all $0\leq t<\theta_0$ and for every $n\in\N$ there exists a $0\leq t<\theta_0$ such that $D\sF_t^n$ has compact values. 
	\end{itemize}
	If $\phi^0$ is a weakly hyperbolic $\theta_1$-periodic solution to $(\Delta_0)$ and there exists a function $\Gamma_0:\R_+\to\R_+$ with $\lim_{\varrho\searrow 0}\Gamma_0(\varrho)=0$ satisfying for all $0\leq t<\theta$ that
	\begin{align}
		\norm{\sF_t^n(\phi_t^0)-\sF_t^0(\phi_t^0)}_X
		&\leq
		\Gamma_0(\tfrac{1}{n}),
		\label{thmsol3}\\
		\lim_{n\to\infty}\norm{D\sF_t^n(\phi_t^0)-D\sF_t^0(\phi_t^0)}_{L(X)}
		&=
		0,
		\label{thmsol1}
	\end{align}
	then there exist reals $\rho_0>0$ and $N_0\in\N$ such that the following hold for all $n\geq N_0$: 
	\begin{enumerate}
		\item There is a unique $\theta$-periodic solution $\phi^n$ to \eqref{deq} in the neighborhood $\cB_{\rho_0}(\phi^0)$, it is weakly hyperbolic and there exists a constant $K_0\geq 0$ such that
		\begin{equation}
			\sup_{t\in\Z}\norm{\phi_t^n-\phi_t^0}_X
			\leq
			K_0\Gamma_0(\tfrac{1}{n}),
			\label{thmsol4}
		\end{equation}

		\item with $\phi^0$ also the solution $\phi^n$ to \eqref{deq} is hyperbolic with the same Morse index. 
	\end{enumerate}
\end{theorem}
As the subsequent proof and \lref{lemhatdera} reveal, the constant $K_0\geq 0$ essentially depends on the distance of the Floquet spectrum of $\phi^0$ to the point $1\in\C$. The value of $K_0$ blows up as this distance shrinks to $0$, i.e.\ when (weak) hyperbolicity is lost. 
\begin{proof}
	Let $I$ denote the identity mapping on the Cartesian product $X^\theta$. Our aim is to apply the Implicit Function \tref{thmifta} with the open set $\varOmega=\hat U$, Banach spaces $\sX=\sY=X^\theta$, the parameter space $\Lambda:=\set{0}\cup\set{\tfrac{1}{n}:\,n\in\N}\subseteq\R$ with metric $d(\lambda_1,\lambda_2):=\abs{\lambda_1-\lambda_2}$, $\lambda_0:=0$, $x_0:=\hat\phi^0$, $y_0:=0$ and the mapping
	\begin{align*}
		T:\hat U\tm\Lambda&\to X^\theta,&
		T(x,\lambda)&:=
		\begin{cases}
			\hat\sF^n(\hat\phi)-\hat\phi,&\lambda=\tfrac{1}{n},\\
			\hat\sF^0(\hat\phi)-\hat\phi,&\lambda=0
		\end{cases}
	\end{align*}
	with $x=\hat\phi$. Let us first verify the assumptions of \tref{thmifta}. It follows from \eqref{fhatdef} that the mapping $T$ is well-defined. 

	\underline{ad (i')}: Thanks to \lref{lemlinb}(a), for the $\theta$-periodic solution $\phi^0$ of $(\Delta_0)$ the resulting tuple $\hat\phi^0$ is a fixed point of $\hat\sF^0$ and therefore $T(x_0,\lambda_0)=\hat\sF^0(\hat\phi^0)-\hat\phi^0=0$. 

	\underline{ad (ii')}: Referring to \lref{lemhatder} and assumption (i) every mapping $\hat\sF^n$ is continuously differentiable and so is each $T(\cdot,\lambda)$, $\lambda\in\Lambda$. Moreover, the partial derivative $D_1T(x_0,\lambda_0)=D\hat\sF^0(\hat\phi^0)-I$ is invertible, because otherwise $1\in\sigma(D\hat\sF^0(\hat\phi^0))$ and thus $1\in\sigma(D\hat\sF^0(\hat\phi^0))^\theta\setminus\set{0}=\sigma(\Xi_\theta^0)\setminus\set{0}$ by \lref{lemhatdera}. This contradicts the weak hyperbolicity assumption on the solution $\phi^0$. 

	\underline{ad (iii')}: First, we obtain \eqref{thmift1} from the estimates
	$$
		\norm{T(x_0,\lambda)-T(x_0,\lambda_0)}
		=
		\norm{\hat\sF^n(\hat\phi^0)-\hat\sF^0(\hat\phi^0)}
		\stackrel{\eqref{fhatdef}}{=}
		\max_{t=1}^\theta\norm{\sF_t^n(\phi_t^0)-\sF_t^0(\phi_t^0)}
	$$
	and thus $\norm{T(x_0,\lambda)-T(x_0,\lambda_0)}\leq\Gamma_0(\lambda)$ (cf.~\eqref{thmsol3}) for all $\lambda=\tfrac{1}{n}\in\Lambda$. Second, by assumption (ii) the derivatives $D\sF_t^n:U_t\to L(X)$ are uniformly continuous on bounded sets, uniformly in $n\in\N$, and consequently there exist moduli of continuity $\omega_t:\R_+\to\R_+$ satisfying $\lim_{\varrho\searrow 0}\omega_t(\varrho)=0$ and
	$$
		\norm{D\sF_t^n(\phi_t)-D\sF_t^n(\phi_t^0)}
		\leq
		\omega_t(\norm{\phi_t-\phi_t^0})\fall n\in\N,\,1\leq t\leq\theta,
	$$
	where $\hat\phi\in\hat U$. By the triangle inequality this results in
	\begin{eqnarray*}
		&&
		\norm{D_1T(x,\lambda)-D_1T(x_0,\lambda_0)}\\
		& \leq &
		\norm{D\hat\sF^n(\hat\phi)-D\hat\sF^n(\hat\phi^0)}+
		\norm{D\hat\sF^n(\hat\phi^0)-D\hat\sF^0(\hat\phi^0)}\\
		& \stackrel{\eqref{fhatdef}}{=} &
		\max_{t=1}^\theta\norm{D\sF_t^n(\phi_t)-D\sF_t^n(\phi_t^0)}+
		\max_{t=1}^\theta\norm{D\sF_t^n(\phi_t^0)-D\sF_t^0(\phi_t^0)}\\
		& \leq &
		\max_{t=1}^\theta\omega_t(\norm{\phi_t-\phi_t^0})+
		\max_{t=1}^\theta\norm{D\sF_t^n(\phi_t^0)-D\sF_t^0(\phi_t^0)}
	\end{eqnarray*}
	for all $\lambda=\tfrac{1}{n}\in\Lambda$. Now, with $\Omega'(\varrho):=\max_{t=1}^\theta\bigl\|D\sF_t^{[1/\varrho]}(\phi_t^0)-D\sF_t^0(\phi_t^0)\bigr\|$ satisfying $\lim_{\varrho\searrow 0}\Omega'(\varrho)=0$ due to \eqref{thmsol1}, this gives for all $\lambda=\tfrac{1}{n}\in\Lambda$ that
	$$
		\norm{D_1T(x,\lambda)-D_1T(x_0,\lambda_0)}
		\stackrel{\eqref{thmsol3}}{\leq}
		\max_{t=1}^\theta\omega_t(\norm{\phi_t-\phi_t^0})+\Omega'(\tfrac{1}{n})
		\leq
		\Gamma\bigl(\bigl\|\hat\phi-\hat\phi^0\bigr\|,\lambda\bigr),
	$$
	with the function $\Gamma(\varrho_1,\varrho_2):=\max_{t=1}^\theta\omega_t(\varrho_1)+\Omega'(\varrho_2)$, which clearly satisfies the limit relation $\lim_{\varrho_1,\varrho_2\searrow 0}\Gamma(\varrho_1,\varrho_2)=0$, i.e.\ \eqref{thmift2} holds.

	(a) Because the assumptions (i'--iii') of \tref{thmifta} hold, we can choose $\rho,\delta>0$ so small that \eqref{thmift3} holds for e.g.\ $q:=\tfrac{1}{2}$. Moreover, there exists a unique fixed point function $\hat\phi:B_\delta(\lambda_0)\to\bar B_{\rho_0}(\hat\phi^0,X^\theta)$ with $\hat\sF^n(\hat\phi(\tfrac{1}{n}))=\hat\phi(\tfrac{1}{n})$ for all $n>\tfrac{1}{\delta}$. Then \lref{lemlinb}(b) guarantees that $\phi^n:=(\overline{\phi_0(\tfrac{1}{n}),\ldots,\phi_{\theta-1}(\tfrac{1}{n})})$ is the desired $\theta$-periodic solution to \eqref{deq} whenever $n\geq N_0:=[\tfrac{1}{\delta}]+1$. We establish that the solutions $\phi^n$ are weakly hyperbolic. For this purpose, let $\eps>0$. First, thanks to \eqref{thmsol1} there exists a $n_1\in\N$ such that
	$$
		\norm{D\sF_t^0(\phi_t^0)-D\sF_t^n(\phi_t^0)}
		\leq
		\tfrac{\eps}{3}\fall t\in\Z,\, n\geq n_1.
	$$
	We know from \tref{thmifta}(c) that $\lim_{n\to\infty}\sup_{t\in\Z}\norm{\phi_t^n-\phi_t^0}=0$ and since $D\sF_t^n$ is equi\-continuous by assumption (ii), there exists a $n_2\in\N$ such that
	$$
		\norm{D\sF_t^n(\phi_t^0)-D\sF_t^n(\phi_t^n)}
		\leq
		\tfrac{\eps}{3}\fall t\in\Z,\, n\geq n_2.
	$$
	Combining the last two inequalities readily yields $\norm{D\sF_t^0(\phi_t^0)-D\sF_t^n(\phi_t^n)}<\eps$ for all $t\in\Z$ and $n\geq\max\set{n_1,n_2}$, which establishes the limit relation
	\begin{equation}
		\lim_{n\to\infty}\norm{D\sF_t^n(\phi_t^n)-D\sF_t^0(\phi_t^0)}=0
		\fall t\in\Z.
		\label{zerostar}
	\end{equation}
	Second, assumption (ii) implies that the period operators $\Xi_\theta^n$ of \eqref{var}, $n\in\N$, contain a compact factor and hence are compact \cite[p.~417, Thm.~1.2]{lang:93}. Thus, \tref{thmhyper}(a) applies to $\sK_t^n:=D\sF_t^n(\phi_t^n)$, $n\in\N_0$, and shows that $\phi^n$ are weakly hyperbolic. Finally, given $N_0$ and $\hat\phi$ as in (a) one has
	$$
		\norm{\phi_t^n-\phi_t^n}
		\stackrel{\eqref{prodnorm}}{\leq}
		\bigl\|\hat\phi^n-\hat\phi^0\bigr\|
		\leq
		K_0\Gamma_0(\tfrac{1}{n})\fall n\geq N_0,\,t\in\Z
	$$
	with $K_0:=2\bigl\|[D\hat\sF^0(\hat\phi^0)-I]^{-1}\bigr\|$, which concludes the proof of (a). 

	(c) In case the solution $\phi^0$ is hyperbolic, then due to \eqref{zerostar} and the compactness of the period operators $\Xi_\theta^n$ (see above), \tref{thmhyper}(b) applies to $\sK_t^n:=D\sF_t^n(\phi_t^n)$, $n\in\N_0$. It follows that the solutions $\phi^n$ are hyperbolic as well. 
	\qed
\end{proof}

The dynamics of difference equations \eqref{deq} in the vicinity of hyperbolic solutions $\phi^n$ is determined by a saddle-point structure consisting of local stable and unstable manifolds resp.\ fiber bundles \cite[p.~256ff, Sect.~4.6]{poetzsche:10b} (in the periodic case). These sets allow a dynamical characterization and, given some $r_0>0$, we define the \emph{local stable fiber bundle}
$$
	\cW_+^n
	:=
	\set{(\tau,u_\tau)\in\cB_{r_0}(\phi^n):\,\vphi^n(t;\tau,u_\tau)-\phi_t^n\xrightarrow[t\to\infty]{}0}
$$
and the \emph{local unstable fiber bundle}
$$
	\cW_-^n
	:=
	\set{(\tau,u_\tau)\in\cB_{r_0}(\phi^n):
	\begin{array}{l}
		\text{there exists a solution $(\phi_t)_{t\leq\tau}$ of }\eqref{deq}\\
		\text{with $\phi_\tau=u_\tau$ and }\phi_t-\phi_t^n\xrightarrow[t\to-\infty]{}0
	\end{array}}
$$
associate to $\phi^n$. The following result relates the fiber bundles of the perturbed equations \eqref{deq}, $n\in\N$, to that of the initial problem $(\Delta_0)$. It requires that $\set{\sF_t^n}_{n\in\N}$ is \emph{equidifferentiable} in each $u\in U_t$, that is there exists a $D\sF_t^n(u)\in L(X)$ such that
$$
	\lim_{h\to 0}\tfrac{1}{\norm{h}_X}\norm{\sF_t^n(u+h)-\sF_t^n(u)-D\sF_t^n(u)h}_X=0\fall t\in\Z
$$
holds uniformly in $n\in\N$. 
\begin{figure}[ht]
	\includegraphics[scale=0.64]{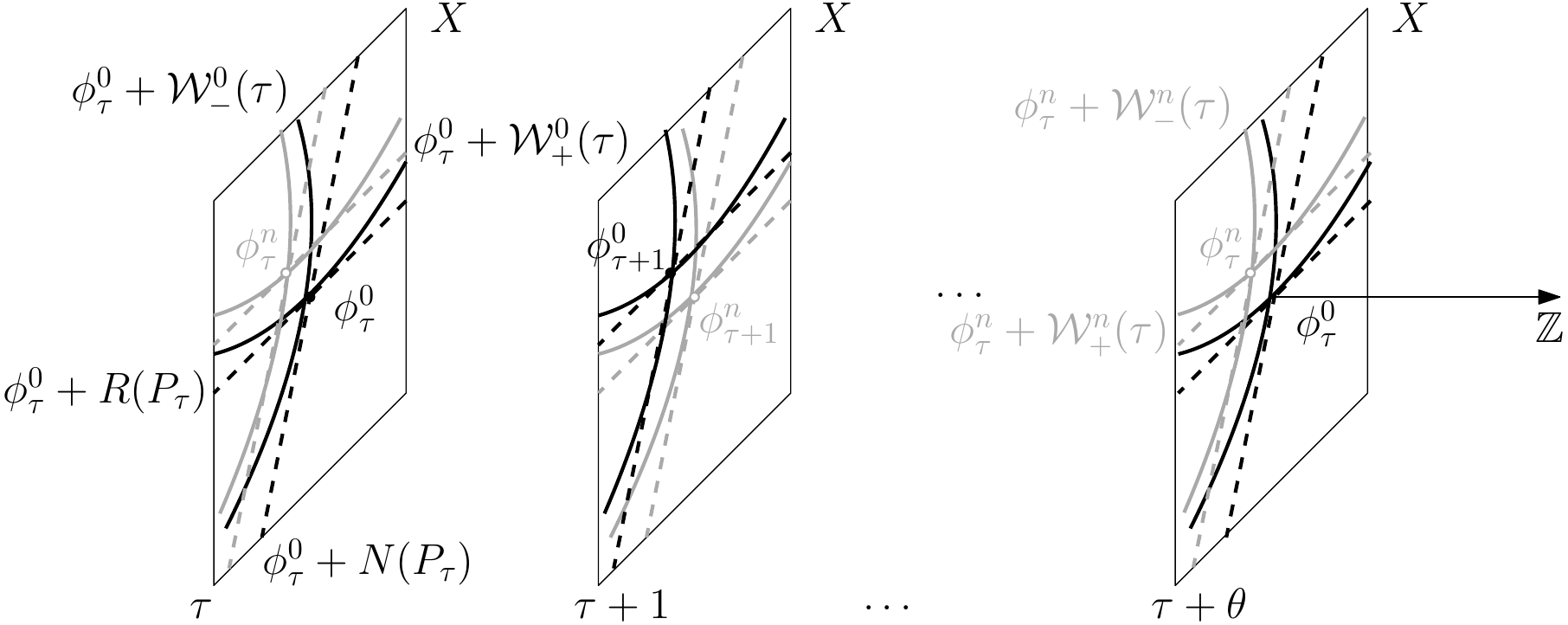}
	\caption{Persistence of the saddle-point structure near a hyperbolic solution: A $\theta_1$-periodic solution $\phi^0$ ($\bullet$ black) of a $\theta_0$-periodic equation $(\Delta_0)$ persists as as $\theta$-periodic hyperbolic solution $\phi^n$ ($\circ$ grey) to \eqref{deq}, $n\geq N_1$ (cf.~\tref{thmsol}). The corresponding stable bundle $\phi^0+\cW_+^0$ (black fibers) persists as $\phi^n+\cW_+^n$ (grey fibers), both are locally graphs over $R(P_\tau)$ (dashed), while the unstable bundle $\phi^0+\cW_-^0$ (black fibers) persists as $\phi^n+\cW_-^n$ (grey fibers), being locally graphs over $N(P_\tau)$ (dashed, cf.~\tref{thmbundle}).}
	\label{figsaddle}
\end{figure}
We can now show that the saddle-point structure near hyperbolic periodic solutions to $(\Delta_0)$ is preserved under perturbation (see \fref{figsaddle}). 
\begin{theorem}[perturbed stable and unstable fiber bundles]\label{thmbundle}
	Let $\theta=\lcm\set{\theta_0,\theta_1}$ and $m\in\N$. Suppose that the $\theta_0$-periodic difference equations \eqref{deq}, $n\in\N_0$, fulfill:
	\begin{itemize}
		\item[(i)] $\sF_t^n:U_t\to X$ are $m$-times continuously differentiable for all $n\in\N_0$ on a convex, open set $U_t$ and $\set{\sF_t^n}_{n\in\N}$ is equidifferentiable for all $0\leq t<\theta_0$, 

		\item[(ii)] $D\sF_t^n:U_t\to L(X)$, $n\in\N$, are uniformly continuous on bounded sets uniformly in $n\in\N$, the family $\set{D\sF_t^n}_{n\in\N}$ is equicontinuous for all $0\leq t<\theta_0$ and for every $n\in\N$ there exists a $0\leq t<\theta_0$ such that $D\sF_t^n$ has compact values. 
	\end{itemize}
	If $\phi^0$ is a hyperbolic $\theta_1$-periodic solution to $(\Delta_0)$ satisfying \eqref{thmsol3}, \eqref{thmsol1} and $(P_t)_{t\in\Z}$ denotes the invariant projector onto the stable vector bundle $\cV_+^0$ of $(V_0)$ (cf.~\tref{thmhyper}), then there exist $\rho_1>0$ and integers $N_1\geq N_0$ so that the following holds for $n\geq N_1$ or $n=0$, and the $\theta$-periodic hyperbolic solutions $\phi^n$ ensured by \tref{thmsol}: 
	\begin{enumerate}
		\item The local stable fiber bundle $\cW_+^n$ of \eqref{deq} allows the representation
		$$
			\cW_+^n
			=
			\phi^n+\set{(\tau,v+w_+^n(\tau,v))\in\Z\tm X:\,v\in B_{\rho_1}(0,R(P_\tau))}
		$$
		as graph of a mapping $w_+^n:\Z\tm X\to X$ with 
		$$
			w_+^n(\tau+\theta,u)=w_+^n(\tau,u)=w_+^n(\tau,P_\tau u)\in N(P_\tau)
			\fall\tau\in\Z
		$$
		and $u\in X$. Moreover, $w_+^n(\tau,0)\equiv 0$ on $\Z$, the Lipschitz mappings $w_+^n(\tau,\cdot)$ are of class $C^m$ and the stable fiber bundles of \eqref{deq} and $(\Delta_0)$ are related via
		\begin{multline}
			\norm{w_+^n(\tau,v)-w_+^0(\tau,v)}_X
			\label{thmbundle3}\\
			\leq
			\frac{4K}{1-\beta}\sup_{\tau\leq t}
			\norm{
			\int_0^1
			\intcc{
			D\sF_t^0(\phi_t^0+\vartheta\phi_t)-D\sF_t^n(\phi_t^n+\vartheta\phi_t)}\phi_t\d\vartheta}_X
		\end{multline}
		for all $\tau\in\Z$, $v\in B_{\rho_1}(0,R(P_\tau))$, where $\phi_t=\vphi^0(t;\tau,\phi_\tau^0+v+w_+^0(\tau,v))-\phi_t^0$ whenever $\tau\leq t$. 

		\item The local unstable fiber bundle $\cW_-^n$ of \eqref{deq} allows the representation
		\begin{equation*}
			\cW_-^n
			=
			\phi^n+\set{(\tau,v+w_-^n(\tau,v))\in\Z\tm X:\,v\in B_{\rho_1}(0,N(P_\tau))}
		\end{equation*}
		as graph of a mapping $w_-^n:\Z\tm X\to X$ with 
		$$
			w_-^n(\tau+\theta,u)=w_-^n(\tau,u)=w_-^n(\tau,[I_X-P_\tau]u)\in R(P_\tau)
			\fall\tau\in\Z
		$$
		and $u\in X$. Moreover, $w_-^n(\tau,0)\equiv 0$ on $\Z$, the Lipschitz mappings $w_-^n(\tau,\cdot)$ are of class $C^m$ and the unstable fiber bundles of \eqref{deq} and $(\Delta_0)$ are related via
		\begin{multline}
			\norm{w_-^n(\tau,v)-w_-^0(\tau,v)}_X
			\label{thmbundle6}\\
			\leq
			\frac{4K}{1-\beta}\sup_{t\leq\tau}
			\norm{
			\int_0^1
			\intcc{
			D\sF_t^0(\phi_t^0+\vartheta\phi_t)-D\sF_t^n(\phi_t^n+\vartheta\phi_t)}\phi_t\d\vartheta}_X
		\end{multline}
		for all $\tau\in\Z$, $v\in B_{\rho_1}(0,N(P_\tau))$, where $(\phi_t)_{t\leq\tau}$ is the (unique) backward solution to $(\Delta_0)$ starting in $(\tau,v+w_-^0(\tau,v))$, and have the same finite dimension. 

		\item $\cW_+^n\cap\cW_-^n=\phi^n$,
	\end{enumerate}
	with the constants $\beta\in(0,1)$, $K\geq 1$ from \tref{thmhyper} applied to $(V_0)$. 
\end{theorem}
In order to achieve convergence as $n\to\infty$ via \eqref{thmbundle3} and \eqref{thmbundle6} one needs the derivatives $D\sF_t^n$ to tend to $D\sF_t^0$ on bounded sets and, thanks to \tref{thmsol}, continuity of the derivative $D\sF_t^0$, $0\leq t<\theta_0$. A concrete illustration follows in Sect.~\ref{sec3}. 
\begin{remark}[alternative representation of $\cW_+^n$ and $\cW_-^n$]
	With some $\tilde\rho_1>0$ the local stable and unstable fiber bundles of $\phi^n$ allow the alternative characterization
	\begin{align*}
		\cW_+^n
		&=
		\phi^n+\set{(\tau,v+\tilde w_+^n(\tau,v))\in\Z\tm X:\,v\in B_{\tilde\rho_1}(0,R(P_\tau^n))},\\
		\cW_-^n
		&=
		\phi^n+\set{(\tau,v+\tilde w_-^n(\tau,v))\in\Z\tm X:\,v\in B_{\tilde\rho_1}(0,N(P_\tau^n))}
	\end{align*}
	as graphs over the vector bundles $\cV_+^n$ resp.\ $\cV_-^n$ of the variational equations \eqref{var}, rather than over the vector bundles $\cV_+^0$ resp.\ $\cV_-^0$ of $(V_0)$ (cf.~\cite[pp.~256ff, Sect.~4.6]{poetzsche:10b}) as in \tref{thmbundle}. In addition, then the associate mappings $\tilde w_+^n(\tau,\cdot),\tilde w_-^n(\tau,\cdot)$ possess values in $N(P_\tau^n)$ resp.\ in $R(P_\tau^n)$ for all $\tau\in\Z$. According to \tref{thmhyper}(c) the corresponding invariant projectors for \eqref{var} satisfy $\lim_{n\to\infty}\norm{P_t^n-P_t^0}=0$ for all $t\in\Z$. Therefore, $\cW_-^n$ and $\cW_-^0$ share their finite dimension. 
\end{remark}
\begin{proof}
	Since the existence of $\cW_+^0,\cW_-^0$ and their properties are well-established in the literature \cite[pp.~187ff]{poetzsche:10b}, we focus on their persistence and the convergence estimates \eqref{thmbundle3} and \eqref{thmbundle6}. Let $\phi^n=(\phi_t^n)_{t\in\Z}$ denote the $\theta$-periodic solutions of \eqref{deq} guaranteed by \tref{thmsol} for $n\geq N_0$. The associate equations of perturbed motion
	\begin{align}
		\tag{$\bar\Delta_n$}
		u_{t+1}&=\bar\sF_t^n(u_t),&
		\bar\sF_t^n(u)&:=\sF_t^n(u+\phi_t^n)-\sF_t^n(\phi_t^n)
		\label{ifb01}
	\end{align}
	are $\theta$-periodic and have the trivial solution. The general solutions $\vphi^n$ of \eqref{deq} and $\bar\vphi^n$ to \eqref{ifb01} are related by $\bar\vphi^n(t;\tau,u)=\vphi^n(t;\tau,u+\phi_\tau^n)-\phi_t^n$ for all $\tau\leq t$. 
	
	(a) For each fixed $\tau\in\Z$ the sequence space
	\begin{equation*}
		\ell_\tau^+
		:=
		\bigl\{(\phi_t)_{\tau\leq t}:\,\phi_t\in X\text{ and }\lim_{t\to\infty}\norm{\phi_t}=0\bigr\},
	\end{equation*}
	is complete w.r.t.\ the $\sup$-norm $\norm{\phi}_\infty:=\sup_{\tau\leq t}\norm{\phi_t}$. For $\bar\rho>0$ so small that $\norm{\phi_t}<\bar\rho$ implies $\phi_t+\phi_t^n\in U_t$ for all $t\in\Z$ and $n\geq N_0$ we introduce the operator 
	\begin{align*}
		T_+^n:B_{\bar\rho}(0,\ell_\tau^+)&\to R(P_\tau)\tm\ell_\tau^+,&
		T_+^n(\phi)_t
		&:=
		\begin{pmatrix}
			P_\tau\phi_\tau,
			\phi_{t+1}-\bar\sF_t^n(\phi_t)
		\end{pmatrix}
	\end{align*}
	for all $\tau\leq t$. Then $u_\tau=P_\tau u_\tau+[I_X-P_\tau]u_\tau\in X$ is contained in the stable bundle of \eqref{ifb01} if and only if $\phi:=\bar\vphi^n(\cdot;\tau,u_\tau)$ satisfies (cf.~\cite[proof of Thm.~3.1]{beyn:87})
	\begin{equation}
		T_+^n(\phi)
		=
		\begin{pmatrix}
			P_\tau u_\tau, 
			0
		\end{pmatrix}.
		\label{ifb03}
	\end{equation}
	Our approach to \eqref{ifb03} using the Lipschitz inverse function \tref{thmlip} is based on the representation $T_+^n=A_++G_+^n$ with
	\begin{align*}
		A_+&\in L(\ell_\tau^+,R(P_\tau)\tm\ell_\tau^+),&
		(A_+\phi)_t&:=
		\begin{pmatrix}
			P_\tau\phi_\tau, 
			\phi_{t+1}-D\sF_t^0(\phi_t^0)\phi_t
		\end{pmatrix},\\
		G_+^n&:\ell_\tau^+\to R(P_\tau)\tm\ell_\tau^+,&
		G_+^n(\phi)_t
		&:=
		\begin{pmatrix}
			0, 
			D\sF_t^0(\phi_t^0)\phi_t-\bar\sF_t^n(\phi_t)
		\end{pmatrix}
	\end{align*}
	for all $\tau\leq t$. Note that the derivatives $D\sF_t^0:U_t\to L(X)$ exist by assumption (i). 

	(I) Claim: \emph{$A_+\in GL(\ell_\tau^+,R(P_\tau)\tm\ell_\tau^+)$ with $\norm{A_+^{-1}}\leq\tfrac{2K}{1-\beta}$.}\\
	First of all, the sequence $(D\sF_t^0(\phi_t^0))_{t\in\Z}$ in $L(X)$ is $\theta$-periodic and therefore $A_+$ is bounded. In order to show that $A_+$ is invertible, given $v_\tau\in R(P_\tau)$ and a sequence $\psi\in\ell_\tau^+$, we observe that $A_+\phi=(v_\tau,\psi)$ has the unique solution
	$$
		\phi_t
		=
		\Phi^0(t,\tau)P_\tau v_\tau+\sum_{s=\tau}^{t-1}\Phi^0(t,s+1)P_s\psi_s
		-
		\sum_{s=t}^\infty\Phi^0(t,s+1)[I_X-P_s]\psi_s
	$$
	in $\ell_\tau^+$ (a proof can be modelled after e.g.\ \cite[pp.~151--152, Thm.~3.5.3(a)]{poetzsche:10b}). Using the dichotomy estimates \eqref{thmhyper2} it is not hard to show
	$
		\norm{\phi_t}
		\leq
		K\norm{v_\tau}+K\tfrac{1+\beta}{1-\beta}\norm{\psi}_\infty
	$
	for all $\tau\leq t$ and therefore $\norm{A_+^{-1}}\leq K+K\tfrac{1+\beta}{1-\beta}=\tfrac{2K}{1-\beta}$. 

	(II) Claim: \emph{There exist $\rho\in(0,\bar\rho]$, $N_1\geq N_0$ such that $\lip G_+^n|_{B_\rho(0)}\leq\tfrac{1-\beta}{4K}$ holds for all $n\geq N_1$.}\\
	Due to the limit relation \eqref{zerostar} in the proof of \tref{thmsol} there is an $N_1\geq N_0$ with
	$$	
		\norm{D\sF_t^0(\phi_t^0)-D\sF_t^n(\phi_t^n)}
		\leq
		\frac{1-\beta}{8K}\fall t\in\Z,\,n\geq N_1. 
	$$
	We next abbreviate $\sH_t^n(u):=D\sF_t^n(\phi_t^n)u-\bar\sF_t^n(u)$. This function is continuously differentiable
	$
		D\sH_t^n(u)
		=
		D\sF_t^n(\phi_t^n)-D\bar\sF_t^n(u)
		=
		D\sF_t^n(\phi_t^n)-D\sF_t^n(u+\phi_t^n).
	$
	The Mean Value Inequality \cite[p.~342, Cor.~4.3]{lang:93} and the fact that $\sF_t^n$ is equidifferentiable by assumption~(i) with continuous derivative thus implies that there exists a $\rho\in(0,\bar\rho]$ such that
	$
		\norm{\sH_t^n(u)-\sH_t^n(\bar u)}
		\leq
		\tfrac{1-\beta}{8K}\norm{u-\bar u}
	$
	for all $t\in\Z$, $u,\bar u\in B_\rho(0,X)$ and $n\geq N_1$. In combination, due to the representation
	$$
		G_+^n(\phi)_t
		=
		\begin{pmatrix}
			0, 
			D\sF_t^0(\phi_t^0)\phi_t-D\sF_t^n(\phi_t^n)\phi_t
		\end{pmatrix}
		+
		\begin{pmatrix}
			0, 
			\sH_t^n(\phi_t)
		\end{pmatrix}
	$$
	we finally obtain for all $\phi,\bar\phi\in B_\rho(0,\ell_\tau^+)$ that
	$$
		\norm{G_+^n(\phi)-G_+^n(\bar\phi)}_\infty
		\leq
		\tfrac{1-\beta}{4K}\norm{\phi-\bar\phi}_\infty\fall n\geq N_1.
	$$

	(III) In this step we apply the Lipschitz inverse function \tref{thmlip} to solve the nonlinear equation \eqref{ifb03} in the Banach spaces $\sX=\ell_\tau^+$, $\sY=R(P_\tau)\tm\ell_\tau^+$, points $x_0:=0$, $y_0:=(P_\tau u_\tau,0)$, the Lipschitz constant $l:=\tfrac{1-\beta}{4K}$ and $\sigma:=\tfrac{1-\beta}{2K}$. Therefore, for every $u_\tau\in B_{\tfrac{1-\beta}{4K^2}\rho}(0,X)$ one has 
	$$
		\norm{(P_\tau u_\tau,0)}
		\stackrel{\eqref{prodnorm}}{=}
		\norm{P_\tau u_\tau}<\frac{1-\beta}{4K}\rho=:\rho_1
	$$
	and there exists a unique solution $\phi_+^n(u_\tau)\in B_\rho(0,\ell_\tau^+)$ to \eqref{ifb03}. Then the function $w_+^n$ parametrizing the stable bundle of \eqref{ifb01} is $w_+^n(\tau,v_\tau):=[I_X-P_\tau]\phi_+^n(v_\tau)_\tau$, where $v_\tau=P_\tau u_\tau$. We define $\phi:=\bar\vphi^0(\cdot;\tau,v_\tau+w_+^0(\tau,v_\tau))$, $\bar\phi:=\phi_+^n(v_\tau)$ and obtain
	\begin{eqnarray*}
		&&
		\norm{w_+^n(\tau,v_\tau)-w_+^0(\tau,v_\tau)}
		=
		\norm{\bar\phi_\tau-\phi_\tau}
		\leq
		\norm{\bar\phi-\phi}_\infty\\
		& \stackrel{\eqref{thmlip1}}{\leq} &
		\frac{4K}{1-\beta}\norm{T_+^n(\bar\phi)-T_+^n(\phi)}_\infty\\
		& \stackrel{\eqref{prodnorm}}{=} &
		\frac{4K}{1-\beta}\max\bigl\{\norm{P_n[\bar\phi_\tau-\phi_\tau]},\sup_{\tau\leq t}\norm{\bar\phi_{t+1}-\bar\sF_t^n(\bar\phi_t)-[\phi_{t+1}-\bar\sF_t^n(\phi_t)]}\bigr\}\\
		& = &
		\frac{4K}{1-\beta}\sup_{\tau\leq t}\norm{\bar\phi_{t+1}-\bar\sF_t^n(\bar\phi_t)-[\phi_{t+1}-\bar\sF_t^n(\phi_t)]}.
	\end{eqnarray*}
	Because $\bar\phi$ solves \eqref{ifb01} and $\phi$ solves $(\bar\Delta_0)$, this simplifies to
	\begin{eqnarray*}
		&&
		\norm{w_+^n(\tau,v_\tau)-w_+^0(\tau,v_\tau)}\\
		& \leq &
		\frac{4K}{1-\beta}\sup_{\tau\leq t}\norm{\phi_{t+1}-\bar\sF_t^n(\phi_t)]}
		=
		\tfrac{4K}{1-\beta}\sup_{\tau\leq t}\norm{\bar\sF_t^0(\phi_t)-\bar\sF_t^n(\phi_t)]}\\
		& \stackrel{\eqref{ifb01}}{=} &
		\frac{4K}{1-\beta}\sup_{\tau\leq t}\norm{\sF_t^0(\phi_t+\phi_t^0)-\sF_t^0(\phi_t^0)-[\sF_t^n(\phi_t+\phi_t^n)+\sF_t^n(\phi_t^n)]}
	\end{eqnarray*}
	and it remains to estimate the right-hand side in this inequality. Since $U_t$ is assumed to be convex, we apply the Mean Value Theorem \cite[p.~341, Thm.~4.2]{lang:93} and arrive at
	\begin{eqnarray*}
		&&
		\norm{w_+^n(\tau,v_\tau)-w_+^0(\tau,v_\tau)}\\
		& \leq &
		\frac{4K}{1-\beta}\sup_{\tau\leq t}
		\norm{\int_0^1[D\sF_t^0(\phi_t^0+\vartheta\phi_t)-D\sF_t^n(\phi_t^n+\vartheta\phi_t)]\phi_t\d\vartheta}
	\end{eqnarray*}
	for all $v_\tau\in R(P_\tau)$. Here, for $\tau\leq t$ one has the relation
	$
		\phi_t
		=
		\bar\vphi^0(t;\tau,v_\tau+w_+^0(\tau,v_\tau))
		=
		\vphi^0(t;\tau,\phi_\tau^0+v_\tau+w_+^0(\tau,v_\tau))-\phi_t^0.
	$

	(b)	The argument is dual to the proof of (a), but now one works in the sequence space
	$
		\ell_\tau^-:=\bigl\{(\phi_t)_{t\leq\tau}:\,\phi_t\in X\text{ and }\lim_{t\to-\infty}\norm{\phi_t}=0\bigr\}
	$
	being complete in the $\sup$-norm. One applies \tref{thmlip} with $\sX=\ell_\tau^-$, $\sY=N(P_\tau)\tm\ell_\tau^-$ and $x_0:=0$, $y_0:=(u_\tau-P_\tau u_\tau,0)$, $l:=\tfrac{1-\beta}{4K}$, $\sigma:=\tfrac{1-\beta}{2K}$ to the nonlinear operator 
	\begin{align*}
		T_-^n:B_\rho(0,\ell_\tau^-)&\to N(P_\tau)\tm\ell_\tau^-,&
		T_-^n(\phi)_t
		&:=
		\begin{pmatrix}
			\phi_\tau-P_\tau\phi_\tau,
			\phi_t-\bar\sF_{t-1}^n(\phi_{t-1})
		\end{pmatrix}
	\end{align*}
	for all $t\leq\tau$. If the unique solution to $T_-^n(\phi)=(P_\tau u_\tau,0)$ is denoted by $\phi_-^n(v_\tau)\in\ell_\tau^-$ such that $v_\tau=u_\tau-P_\tau u_\tau$, then $w_-^n(\tau,v_\tau):=P_\tau\phi_-^n(v_\tau)_\tau$ has the claimed properties. 
	
	(c) is a consequence of \cite[pp.~259--260, Thm.~4.6.4]{poetzsche:10b}. 
	\qed
\end{proof}
\section{Urysohn integrodifference equations}
\label{sec3}
Let us now illustrate the applicability of our abstract perturbation results from Sect.~\ref{sec2}, when the initial problem $(\Delta_0)$ is an integrodifference equation 
\begin{align}
	\tag{$I_0$}
	u_{t+1}&=\sF_t^0(u_t),&
	\sF_t^0(u)&:=\int_\Omega f_t(\cdot,y,u(y))\d y, 
\end{align}
whose right-hand side is an Urysohn operator over a compact nonempty $\Omega\subset\R^\kappa$. For the sake of having well-defined and smooth mappings $\sF_t^0$, $t\in\Z$, in an ambient setting, several assumptions on the \emph{kernel functions} $f_t$ are due: 
\begin{hypothesis}\label{hyp}
	Let $m\in\N$ and $\alpha\in(0,1]$. Suppose there exists a $\theta_0\in\N$ and open, convex sets $Z_t\subseteq\R^d$ such that the kernel functions
	\begin{align}
		f_t&=f_{t+\theta_0}:\Omega^2\tm \overline{Z_t}\to\R^d,&
		Z_t&=Z_{t+\theta_0}\fall t\in\Z
		\label{hyp00}
	\end{align}
	fulfill the following assumptions for all $0\leq t<\theta_0$ and $0\leq k\leq m$:
	\begin{itemize}
		\item[(i)] The derivative $D_3^kf_t:\Omega^2\tm \overline{Z_t}\to L_k(\R^d)$ exists as continuous function, 

		\item[(ii)] for all $r>0$ there exists a continuous function $h_r:\Omega\to\R_+$ such that
		\begin{align}
			\abs{D_3^kf_t(x,y,z)-D_3^kf_t(\bar x,y,z)}_{L_k(\R^d)}&\leq h_r(y)\abs{x-\bar x}^\alpha
			\label{hyp01}
		\end{align}
		for all $x,\bar x,y\in\Omega$, $z\in \overline{Z_t}\cap\bar B_r(0)$, 

		\item[(iii)] for all $r>0$ there exists a function $c_r:\R_+\tm\Omega\to\R_+$ satisfying the limit relation $\lim_{\delta\searrow 0}\sup_{y\in\Omega}c_r(\delta,y)=0$, such that $\abs{z-\bar z}\leq\delta$ implies
		\begin{multline}
			\abs{D_3^kf_t(x,y,z)-D_3^kf_t(x,y,\bar z)-\intcc{D_3^kf_t(\bar x,y,z)-D_3^kf_t(\bar x,y,\bar z)}}_{L_k(\R^d)}\\
			\leq
			c_r(\delta,y)\abs{x-\bar x}^\alpha
			\fall x,\bar x,y\in\Omega,\,\bar z\in\overline{Z_t}\cap\bar B_r(0).
			\label{hyp03}
		\end{multline}
	\end{itemize}
\end{hypothesis}
Let $C(\Omega,\R^d)$ denote the set of continuous functions $u:\Omega\to\R^d$ equipped with the norm $\norm{u}_0:=\sup_{x\in\Omega}\abs{u(x)}$. If $\alpha\in(0,1]$, then functions $u:\Omega\to\R^d$ having a bounded H\"older constant
$$
	[u]_\alpha:=\sup_{\substack{x,\bar x\in\Omega\\ x\neq \bar x}}\frac{\abs{u(x)-u(\bar x)}}{\abs{x-\bar x}^\alpha}<\infty
$$
are called $\alpha$-\emph{H\"older} (\emph{Lipschitz} in case $\alpha=1$) and $C^\alpha(\Omega,\R^d)\subset C(\Omega,\R^d)$ denotes the entity of all such functions. It is a Banach space when equipped with the norm
\begin{equation}
	\norm{u}_\alpha
	:=
	\begin{cases}
		\norm{u}_0,&\alpha=0,\\
		\max\set{\norm{u}_0,[u]_\alpha},&\alpha\in(0,1].
	\end{cases}
	\label{normdef}
\end{equation}
Since the compact domain $\Omega$ is fixed throughout, we conveniently abbreviate
\begin{align*}
	\Cda&:=C^\alpha(\Omega,\R^d),&
	\Cd&:=C(\Omega,\R^d)
\end{align*}
and obtain the open sets $U_t:=\set{u\in\Cda:\,u(\Omega)\subset Z_t}$ for all $t\in\Z$. 

For our subsequent analysis it is important to note that \href{hyp} implies the corresponding assumptions made in \cite[Sect.~2]{poetzsche:20b}. In detail, one has:
\begin{proposition}[properties of \eqref{ide0}]\label{propU0}
	Let $t\in\Z$. If \href{hyp} holds, then the Urysohn operator $\sF_t^0=\sF_{t+\theta_0}^0:U_t\to\Cda$ is well-defined, completely continuous and of class $C^m$ with compact derivative 
	\begin{equation}
		D\sF_t^0(u)v
		=
		\int_\Omega D_3f_t(\cdot,y,u(y))v(y)\d y
		\fall u\in U_t,\,v\in\Cda.
		\label{propDU0a}
	\end{equation}
\end{proposition}
Combined with the solution identity this shows that entire solutions $\phi$ to \eqref{ide0} inherit the smoothness of the kernel function, i.e.\ $\phi_t\in\Cda$, $t\in\Z$. Yet for kernel functions of convolution type a higher smoothness can be expected (cf.~\cite[Sect.~2.3]{poetzsche:20b}). 
\begin{proof}
	Above all, \eqref{ide0} and \eqref{hyp00} show that $\sF_t^0$ is $\theta_0$-periodic in $t$. The results from \cite{poetzsche:20b} formulated in an abstract measure-theoretical set-up apply to $\sF_t^0$ with the $\kappa$-dimen\-sional Lebesgue measure $\mu=\lambda_\kappa$. By \cite[Thm.~2.6]{poetzsche:20b}, $\sF_t^0$ is well-defined and due to \cite[Cor.~2.7(i)]{poetzsche:20b} also completely continuous. In \cite[Thm.~2.12]{poetzsche:20b} it is shown that $\sF_t^0$ is of class $C^m$ and \cite[p.~89, Prop.~6.5]{martin:76} implies that $D\sF_t^0(u)$, $u\in U_t$, is compact. 
	\qed
\end{proof}
\begin{corollary}\label{corlip}
	Let $t\in\Z$ and $2\leq m$. If for every $r>0$ there exists a continuous function $l_r:\Omega^2\to\R_+$ with
	$$
		\abs{D_3f_t(x,y,z)-D_3f_t(x,y,\bar z)}_{L(\R^d)}\leq l_r(x,y)\abs{z-\bar z}
		\fall x,y\in\Omega
	$$
	and $z,\bar z\in Z_t\cap\bar B_r(0)$, then $D\sF_t^0:U_t\to L(\Cda)$ is Lipschitz on $C_d^0$-bounded sets, that is, for each $r>0$ there exists a $L_r\geq 0$ such that
	\begin{equation}
		\norm{D\sF_t^0(u)-D\sF_t^0(\bar u)}_{L(\Cda)}
		\leq
		L_r\norm{u-\bar u}_\alpha
		\fall u,\bar u\in U_t\cap\bar B_r(0,\Cd)
		\label{corlip2}
	\end{equation}
	with the Lipschitz constant $L_r:=\max\set{\sup_{\xi\in\Omega}\int_{\Omega} l_r(\xi,y)\d y,\int_\Omega h_r(y)\d y}$. 
\end{corollary}
\begin{proof}
	Let $v\in\Cda$, $r>0$ and $u,\bar u\in U_t\cap\bar B_r(0,\Cd)$. 
	
	(I) We derive that
	\begin{align*}
		\abs{[D\sF_t^0(u)-D\sF_t^0(\bar u)]v(x)}
		&\stackrel{\eqref{propDU0a}}{\leq}
		\int_{\Omega}\abs{D_3f_t(x,y,u(y))-D_3f_t(x,y,\bar u(y))}\abs{v(y)}\d y\\
		&\leq
		\int_{\Omega}l_r(x,y)\abs{u(y)-\bar u(y)}\abs{v(y)}\d y\\
		&\stackrel{\eqref{normdef}}{\leq}
		\sup_{\xi\in\Omega}\int_{\Omega} l_r(\xi,y)\d y\norm{u-\bar u}_\alpha\norm{v}_\alpha
		\fall x\in\Omega. 
	\end{align*}
	Thus, $\norm{[D\sF_t^0(u)-D\sF_t^0(\bar u)]v}_0\leq\sup_{\xi\in\Omega}\int_\Omega l_r(\xi,y)\d y\norm{u-\bar u}_\alpha\norm{v}_\alpha$ after passing to the least upper bound over all $x\in\Omega$. 
	
	(II) With $Z_t\subseteq\R^d$ also $U_t\subseteq\Cda$ is convex. Therefore, the Mean Value Theorem \cite[p.~341, Thm.~4.2]{lang:93} applies and shows for $x,\bar x\in\Omega$ that
	\begin{align*}
		&
		[D\sF_t^0(u)-D\sF_t^0(\bar u)]v(x)-\intcc{D\sF_t^0(u)-D\sF_t^0(\bar u)}v(\bar x)\\
		\stackrel{\eqref{propDU0a}}{=}&
		\int_{\Omega}\bigl[D_3f_t(x,y,u(y))-D_3f_t(x,y,\bar u(y))\\
		&\quad-\bigl(D_3f_t(\bar x,y,u(y))-D_3f_t(\bar x,y,\bar u(y))\bigr)\bigl]v(y)\d y\\
		=&
		\int_{\Omega}
		\int_0^1\bigl[
		D_3^2f_t(x,y,\bar u(y)+\vartheta(u(y)-\bar u(y)))\\
		&\quad-D_3^2f_t(\bar x,y,\bar u(y)+\vartheta(u(y)-\bar u(y)))\d\vartheta\intcc{u(y)-\bar u(y)}\bigr]v(y)\d y.
	\end{align*}
	Consequently \href{hyp}(ii) leads to
	\begin{align*}
		&
		\abs{[D\sF_t^0(u)-D\sF_t^0(\bar u)]v(x)-\intcc{D\sF_t^0(u)-D\sF_t^0(\bar u)}v(\bar x)}\\
		\stackrel{\eqref{normdef}}{\leq}&
		\int_{\Omega}
		\int_0^1
		\bigl|D_3^2f\bigl(x,y,\bar u(y)+\vartheta(u(y)-\bar u(y))\bigr)\\
		&\quad-D_3^2f\bigl(\bar x,y,\bar u(y)+\vartheta(u(y)-\bar u(y))\bigr)\bigr|\d\vartheta\d y
		\norm{u-\bar u}_\alpha\norm{v}_\alpha\\
		\stackrel{\eqref{hyp01}}{\leq}&
		\int_{\Omega}h_r(y)\d y\norm{u-\bar u}_\alpha\norm{v}_\alpha\abs{x-\bar x}^\alpha
		\fall x,\bar x\in\Omega, 
	\end{align*}
	which guarantees $[[D\sF_t^0(u)-D\sF_t^0(\bar u)]v]_\alpha\leq\int_{\Omega}h_r(y)\d y\norm{u-\bar u}_\alpha\norm{v}_\alpha$. Referring to \eqref{normdef} this implies the local Lipschitz estimate \eqref{corlip2}. 
	\qed
\end{proof}

Along with IDEs \eqref{ide0} we now consider their Nystr\"om discretizations. They are based on \emph{quadrature} ($\kappa=1$) or \emph{cubature rules} ($\kappa>1$), i.e.\ a family of mappings
\begin{align}
	\tag{$Q_n$}
	Q^n:\Cd&\to\R^d,&
	Q^nu&:=\sum_{\eta\in\Omega_n}w_\eta u(\eta)\fall n\in\N
	\label{quad}
\end{align}
determined by a \emph{grid} $\Omega_n\subset\Omega$ of finitely many \emph{nodes} $\eta\in\Omega_n$ and \emph{weights} $w_\eta\geq 0$; the dependence of $w_\eta$ on $n\in\N$ is suppressed here. A rule \eqref{quad} is called (cf.~\cite{hackbusch:95})
\begin{itemize}
	\item \emph{convergent}, if $\lim_{n\to\infty}Q^nu=\int_\Omega u(y)\d y$ holds for all $u\in\Cd$, 

	\item \emph{stable}, provided the weights satisfy
	\begin{align}
		W&:=\sup_{n\in\N}W_n<\infty,&
		W_n&:=\sum_{\eta\in\Omega_n}w_\eta. 
		\label{noa}
	\end{align}
\end{itemize}
Thanks to \cite[p.~20, Thm.~1.4.17]{hackbusch:95}, convergence implies stability. 

In order to evaluate the right-hand side of \eqref{ide0} approximately, we replace the integral by a convergent integration rule \eqref{quad}, $n\in\N$. The resulting \emph{Nystr\"om method} (see \cite{atkinson:92,hackbusch:95} for integral equations) yields the family of difference equations
\begin{align}
	\tag{$I_n$}
	u_{t+1}&=\sF_t^n(u_t),&
	\sF_t^n(u)&:=\sum_{\eta\in\Omega_n}w_\eta f_t(\cdot,\eta,u(\eta)). 
\end{align}
\begin{proposition}[properties of \eqref{iden}]\label{propUn}
	Let $t\in\Z$. If \href{hyp} holds, then the discrete Urysohn operator $\sF_t^n=\sF_{t+\theta_0}^n:U_t\to\Cda$, $n\in\N$, is well-defined, completely continuous and of class $C^m$ with compact derivative 
	\begin{equation}
		D\sF_t^n(u)v
		=
		\sum_{\eta\in\Omega_n}w_\eta D_3f_t(\cdot,\eta,u(\eta))v(\eta)
		\fall u\in U_t,\,v\in\Cda.
		\label{propUna}
	\end{equation}
	Moreover, if \eqref{quad} is stable, then $\set{\sF_t^n}_{n\in\N}$ is equidifferentiable, $D\sF_t^n$ are uniformly continuous on bounded sets uniformly in $n\in\N$ and $\set{D\sF_t^n}_{n\in\N}$ is equicontinuous. 
\end{proposition}
\begin{proof}
	The grids $\Omega_n$, $n\in\N$, are a family of compact and discrete subsets of $\Omega$. If we equip them with the measure $\mu(\Omega_n):=\sum_{\eta\in\Omega_n}w_\eta$, then due to \cite[Ex.~2.2 and Rem.~2.5]{poetzsche:20b} the abstract measure-theoretical integral from \cite{poetzsche:20b} becomes
	$$
		\int_{\Omega_n}f_t(x,y,u(y))\d\mu(y)
		=
		\sum_{\eta\in\Omega_n}w_\eta f_t(x,\eta,u(\eta))\fall x\in\Omega
	$$
	and leads to the discrete integral operators in \eqref{iden}. Given this, well-definedness, complete continuity and smoothness of $\sF_t^n$ result from \cite{poetzsche:20b} as in the proof of \pref{propU0}. From now on, assume that \eqref{quad} is stable and choose $u\in U_t$. 

	(I) \emph{Claim}: $\set{\sF_t^n}_{n\in\N}$ is equidifferentiable.\\
	For functions $h\in\Cda$ the remainder terms \cite[(16) resp.~(18)]{poetzsche:20b} become
	\begin{align*}
		r_0(h)&=\sup_{\vartheta\in[0,1]}\bigg\|\sum_{\eta\in\Omega_n}w_\eta\intcc{D_3f_t(\cdot,\eta,(u+\vartheta h)(\eta))-D_3f_t(\cdot,\eta,u(\eta))}\bigg\|_0,\\
		\rho_0(h)&=\int_0^1\sum_{\eta\in\Omega_n}w_\eta\bar c_r^1(\vartheta\norm{h}_0,y)\d\vartheta
		\leq
		\sum_{\eta\in\Omega_n}w_\eta\bar c_r^1(\norm{h}_0,\eta).
	\end{align*}
	Now it follows from \eqref{noa} that $\lim_{h\to\infty}r_0(h)=\lim_{h\to\infty}\rho_0(h)=0$ hold uniformly in $n\in\N$. This yields the claimed equidifferentiability. 
	
	(II) \emph{Claim}: $D\sF_t^n$ are uniformly continuous on bounded sets uniformly in $n\in\N$ (and thus $\set{D\sF_t^n}_{n\in\N}$ is equicontinuous).\\
	Let $\eps>0$, $v\in\Cda$ and given $u,\bar u\in U_t$ choose $r>0$ so large that $\norm{u}_0,\norm{\bar u}_0\leq r$ holds. Because the (extended) derivative $D_3f_t:\Omega^2\tm\overline{Z_t}\to L(\R^d)$ is uniformly continuous on the compact set $\Omega^2\tm(\overline{Z_t}\cap\bar B_r(0))$, there exists a $\delta_1>0$ such that 
	$$
		\abs{z-\bar z}<\delta_1
		\Rightarrow
		\abs{D_3f_t(x,y,z)-D_3f_t(x,y,\bar z)}<\tfrac{\eps}{2W}
		\fall z,\bar z\in Z_t\cap\bar B_r(0)
	$$
	and $x,y\in\Omega$. If $u,\bar u\in U_t$ satisfy $\norm{u-\bar u}_0<\delta_1$, then we obtain $\abs{u(y)-\bar u(y)}<\delta_1$ for all $y\in\Omega$. First, this implies
	\begin{eqnarray*}
		&&
		\abs{[D\sF_t^n(u)-D\sF_t^n(\bar u)]v(x)}\\
		& \stackrel{\eqref{propUna}}{\leq} &
		\sum_{\eta\in\Omega_n}w_\eta
		\abs{D_3f_t(x,\eta,u(\eta))-D_3f_t(x,\eta,\bar u(\eta))}
		\abs{v(\eta)}\\
		& \stackrel{\eqref{normdef}}{\leq} &
		\sum_{\eta\in\Omega_n}w_\eta\tfrac{\eps}{2W}\norm{v}_\alpha
		\leq
		\tfrac{\eps}{2}\norm{v}_\alpha\fall x\in\Omega
	\end{eqnarray*}
	and passing to the supremum over $x\in\Omega$ yields
	$
		\norm{[D\sF_t(u)-D\sF_t(\bar u)]v}_0\leq\tfrac{\eps}{2}\norm{v}_\alpha.
	$
	Second, from \href{hyp}(iii) there exists a $\delta_2>0$ such that $\sup_{y\in\Omega}c_r(\delta,y)<\tfrac{\eps}{2W}$ for every $\delta\in(0,\delta_2]$ and consequently $\norm{u-\bar u}_0<\delta_2$ guarantees for all $x,\bar x\in\Omega$ that
	\begin{eqnarray*}
		&&
		\abs{[D\sF_t^n(u)-D\sF_t^n(\bar u)]v(x)-[D\sF_t^n(u)-D\sF_t^n(\bar u)]v(\bar x)}\\
		& \stackrel{\eqref{propUna}}{\leq} &
		\sum_{\eta\in\Omega_n}w_\eta
		|D_3f_t(x,\eta,u(\eta))-D_3f_t(x,\eta,\bar u(\eta))\\
		&&
		\qquad
		-[D_3f_t(\bar x,\eta,u(\eta))-D_3f_t(\bar x,\eta,\bar u(\eta))]|
		\abs{v(\eta)}\\
		& \stackrel{\eqref{hyp03}}{\leq} &
		\sum_{\eta\in\Omega_n}w_\eta c_r(\delta,\eta)\abs{x-\bar x}^\alpha\norm{v}_\alpha
		\stackrel{\eqref{noa}}{\leq}
		W\sup_{y\in\Omega}c_r(\delta,\eta)\abs{x-\bar x}^\alpha\norm{v}_\alpha
	\end{eqnarray*}
	and therefore $[[D\sF_t^n(u)-D\sF_t^n(\bar u)]v]_\alpha\leq\tfrac{\eps}{2}\norm{v}_\alpha$. Referring to \eqref{normdef} this results in
	$$
		\norm{u-\bar u}_0<\min\set{\delta_1,\delta_2}
		\quad\Rightarrow\quad
		\norm{[D\sF_t^n(u)-D\sF_t^n(\bar u)]v}_\alpha\leq\tfrac{\eps}{2}\norm{v}_\alpha
	$$
	for all $n\in\N$. Since $v\in\Cda$ was arbitrary, this readily implies the claim. 
	\qed
\end{proof}
\subsection{H\"older continuous kernel functions}
We say an integration rule \eqref{quad} has \emph{consistency order } $\alpha\in(0,1]$ (cf.~\cite[p.~21, Def.~1.4.19]{hackbusch:95}), if there exists a $c_0\geq 0$ with
$$
	\abs{\int_\Omega u(y)\d y-Q^nu}\leq\frac{c_0}{n^\alpha}\norm{u}_\alpha
	\fall u\in\Cda. 
$$
\begin{example}[quadrature rules]
	Let $\Omega=[a,b]$ and $n\in\N$. 
	The \emph{(left)} resp.\ \emph{(right) rectangular rules}
	$Q_{LR}^nu:=\tfrac{b-a}{n}\sum_{j=0}^{n-1}u(a+j\tfrac{b-a}{n})$ and 
	$Q_{RR}^nu:=\tfrac{b-a}{n}\sum_{j=1}^nu(a+j\tfrac{b-a}{n})$ are convergent and satisfy the quadrature error (cf.~\cite[p.~52, Theorem]{davis:rabinowitz:84})
	$$
		\abs{\int_a^bu(y)\d y-Q_i^nu}
		\leq
		\frac{(b-a)^{\alpha+1}}{n^\alpha}[u]_\alpha
		\quad\text{for }i\in\set{LR,RR}.
	$$
	Also the \emph{midpoint rule}
	$
		Q_M^nu:=\tfrac{b-a}{n}\sum_{j=0}^{n-1}u(a+(j+\tfrac{1}{2})\tfrac{b-a}{n})
	$
	is convergent and as in \cite[p.~52, Theorem]{davis:rabinowitz:84} one derives the quadrature error 
	$$
		\abs{\int_a^bu(y)\d y-Q_M^nu}
		\leq
		\frac{(b-a)^{\alpha+1}}{2^\alpha n^\alpha}[u]_\alpha. 
	$$
	The \emph{trapezoidal rule} $Q_T^nu:=\tfrac{1}{2}(Q_{LR}^nu+Q_{RR}^nu)$ is convergent with the same quadrature error as for the rectangular rules. Finally, let $n\in\N$ be even. Representing the \emph{Simpson rule} as convex combination $Q_S^nu:=\tfrac{2}{3}Q_M^{n/2}u+\tfrac{1}{3}Q_T^{n/2}u$, one obtains
	$$
		\abs{\int_a^bu(y)\d y-Q_S^nu}
		\leq
		\frac{2+2^\alpha}{3}
		\frac{(b-a)^{\alpha+1}}{n^\alpha}[u]_\alpha. 
	$$
\end{example}

The next two results provide sufficient conditions on the kernel functions $f_t$ such that the assumptions \eqref{thmsol3} or \eqref{thmsol1} are satisfied for Nystr\"om discretizations \eqref{iden}. 
\begin{proposition}[convergence of $\sF_t^n$]\label{lemnconv}
	Let $t\in\Z$. Suppose \href{hyp} holds and that for every $r>0$ there exists a $l_r^0\geq 0$ such that
	\begin{equation}
		\abs{f_t(x,y,z)-f_t(x,\bar y,\bar z)}
		\leq
		l_r^0\max\set{\abs{y-\bar y}^\alpha,\abs{z-\bar z}}
		\label{lemnconv1}
	\end{equation}
	for all $x,y,\bar y\in\Omega$ and $z,\bar z\in Z_t\cap\bar B_r(0)$.
	If \eqref{quad} has consistency order $\alpha$, then for every $r>0$ there exists a $c_r^0\geq 0$ such that
	\begin{equation}
		\norm{\sF_t^n(u)-\sF_t^0(u)}_\alpha
		\leq
		\frac{c_0c_r^0}{n^\alpha}\fall n\in\N,\,u\in U_t\cap\bar B_r(0,\Cda).
		\label{lemnconv2}
	\end{equation}
\end{proposition}
The magnitude of the constant $c_r^0$ is increasing in the H\"older norm of $u\in U_t$. 
\begin{proof}
	Let $t\in\Z$, $r>0$ and $u\in U_t\cap\bar B_r(0,\Cda)$. Because \eqref{quad} has consistency order $\alpha$, there exists a $c_0\geq 0$ such that
	$
		\abs{\sF_t^0(u)(x)-\sF_t^n(u)(x)}
		\leq
		\tfrac{c_0}{n^\alpha}\norm{f_t(x,\cdot,u(\cdot))}_\alpha
	$
	for all $x\in\Omega$. 
	First, one has $\norm{f_t(x,\cdot,u(\cdot))}_0\leq\sup_{\xi,y\in\Omega}\abs{f_t(\xi,y,u(y))}=:b_t$ for every $x\in\Omega$. Second, due to the assumption \eqref{lemnconv1} we conclude
	$$
		\abs{f_t(x,y,u(y))-f_t(x,\bar y,u(\bar y))}
		\leq
		l_r^0\max\set{1,[u]_\alpha}\abs{y-\bar y}^\alpha
		\fall y,\bar y\in\Omega
	$$
	and therefore $[f_t(x,\cdot,u(\cdot))]_\alpha\leq l_r^0\max\set{1,[u]_\alpha}$ holds. In conclusion, because of \eqref{normdef} we arrive at $\norm{f_t(x,\cdot,u(\cdot))}_\alpha\leq\max\set{b_t,l_r^0\max\set{1,[u]_\alpha}}$ for every $x\in\Omega$ and consequently choose $c_r^0:=\max_{t=1}^{\theta_0}\set{b_t,l_r^0\max\set{1,r}}$. 
	\qed
\end{proof}

\begin{proposition}[convergence of $D\sF_t^n$]\label{lemest}
	Let $t\in\Z$. Suppose \href{hyp} holds and that for every $r>0$ there exist constants
	\begin{itemize}
		\item[(iv)] $l_r^1\geq 0$ such that for all $x,y,\bar y\in\Omega$ and $z\in Z_t\cap\bar B_r(0)$ one has
		$$
			\abs{D_3f_t(x,y,z)-D_3f_t(x,\bar y,\bar z)}_{L(\R^d)}
			\leq
			l_r^1\max\set{\abs{y-\bar y}^\alpha,\abs{z-\bar z}}, 
		$$

		\item[(v)] $\gamma_r\geq 0$ such that for all $x,\bar x,y,\bar y\in\Omega$ and $u\in U_t\cap\bar B_r(0,\Cda)$ one has
		\begin{multline}
			\hspace*{-10mm}
			\abs{D_3f_t(x,y,u(y))-D_3f_t(\bar x,y,u(y))-\intcc{D_3f_t(x,\bar y,u(\bar y))-D_3f_t(\bar x,\bar y,u(\bar y))}}_{L(\R^d)}\\
			\leq
			\gamma_r\abs{x-\bar x}^\alpha\abs{y-\bar y}^\alpha.
			\label{lemest2}
		\end{multline}
	\end{itemize}
	If \eqref{quad} has consistency order $\alpha$, then for every $r>0$ there exists a $c_r^1\geq 0$ such that
	\begin{equation}
		\norm{D\sF_t^n(u)-D\sF_t^0(u)}_{L(\Cda)}
		\leq
		\frac{c_0c_r^1}{n^\alpha}
		\fall n\in\N,\,u\in U_t\cap\bar B_r(0,\Cda).
		\label{lemest3}
	\end{equation}
\end{proposition}
Sufficient conditions for \eqref{lemest2} to hold were given in \cite[Rem.~1]{poetzsche:20} on convex $\Omega\subset\R^\kappa$. Furthermore, the explicit form of the constant $c_r^1$ can be obtained from \cite[(11)]{poetzsche:20}. 
\begin{proof}
	Let $t\in\Z$, $r>0$ and $u\in U_t\cap\bar B_r(0,\Cda)$ be fixed. By \pref{propU0} the derivative of $\sF_t^0$ is $D\sF_t^0(u)v=\int_\Omega D_3f_t(\cdot,y,u(y))v(y)\d y$ for all $v\in\Cda$. Given this, our goal is to apply the convergence result \cite[Thm.~2]{poetzsche:20} with the corresponding kernel $k_t(x,y):=D_3f_t(x,y,u(y))$, whose assumptions are verified next: 

	\underline{ad (i)}: Thanks to $\abs{u(y)}\leq r$ it holds 
	$\abs{k_t(x,y)-k_t(\bar x,y)}\leq h_r(y)\abs{x-\bar x}^\alpha$ for all $x,\bar x\in\Omega$ due to \eqref{hyp01} and therefore $[k_t(\cdot,y)]_\alpha\leq\sup_{\eta\in\Omega}h_r(\eta)$ for all $y\in\Omega$. 

	\underline{ad (ii)}: The assumption (iv) and $[u]_\alpha\leq r$ yield
	$$
		\abs{k_t(x,y)-k_t(x,\bar y)}
		\leq
		l_r^1\max\set{1,[u]_\alpha}\abs{y-\bar y}^\alpha
		\fall y,\bar y\in\Omega
	$$
	and thus $[k_t(x,\cdot)]_\alpha\leq l_r^1\max\set{1,r}$ for all $x\in\Omega$. 

	\underline{ad (iii)}: As consequence of our assumption \eqref{lemest2} one obtains for $x,\bar x,y,\bar y\in\Omega$ that $\abs{k_t(x,y)-k_t(\bar x,y)-\intcc{k_t(x,\bar y)-k_t(\bar x,\bar y)}}\leq\gamma_r\abs{x-\bar x}^\alpha\abs{y-\bar y}^\alpha$.
	
	Finally, combining (i--iii) with the consistency order $\alpha$ of \eqref{quad} shows \eqref{lemest3}. 
	\qed
\end{proof}

Combining the assumptions of \pref{propU0}--\ref{lemest} and \cref{corlip} yields
\begin{corollary}[saddle-point structure of \eqref{ide0}, $C^\alpha$-case]\label{cornyst}
	Suppose \href{hyp} holds with $2\leq m$ and $\phi^0$ is a weakly hyperbolic $\theta_1$-periodic solution to \eqref{ide0}. If
	\begin{enumerate}
		\item[(iv)] \eqref{quad} has consistency order $\alpha\in(0,1]$, 
	\end{enumerate}
	then there exist constants $K_\ast,K_+,K_-\geq 0$ and $N_1\in\N$ such that the following holds for all $n\geq N_1$: 
	The associate weakly hyperbolic and $\theta$-periodic solutions $\phi^n$ to \eqref{iden} satisfy
	\begin{equation}
		\sup_{t\in\Z}\norm{\phi_t^n-\phi_t^0}_\alpha
		\leq
		\frac{K_\ast c_0}{n^\alpha}.
		\label{box1}
	\end{equation}
	If $\phi^0$ is even hyperbolic, then for each $\tau\in\Z$ one has the estimates
	\begin{enumerate}
		\item $\norm{w_+^n(\tau,v)-w_+^0(\tau,v)}_\alpha\leq
		\tfrac{4K}{1-\beta}\tfrac{K_+}{n^\alpha}\sup_{\tau\leq t}\norm{\phi_t}_\alpha$ for all $v\in B_{\rho_1}(0,R(P_\tau))$, 

		\item $\norm{w_-^n(\tau,v)-w_-^0(\tau,v)}_\alpha\leq
		\tfrac{4K}{1-\beta}\tfrac{K_-}{n^\alpha}\sup_{t\leq\tau}\norm{\phi_t}_\alpha$ for all $v\in B_{\rho_1}(0,N(P_\tau))$,
	\end{enumerate}
	with the forward resp.\ backward solution $\phi$ to the IDE \eqref{ide0} from \tref{thmbundle}. 
\end{corollary}
\begin{proof}
	Let $t\in\Z$ and $r:=\max_{t=1}^{\theta_1}\norm{\phi_t^0}_\alpha$. It results from \pref{propU0} and \ref{propUn} that \eqref{iden}, $n\in\N_0$, satisfy the assumptions (i), (ii) of \tref{thmsol}. Moreover, \pref{lemnconv} implies \eqref{thmsol3} with $\Gamma_0(\varrho):=c_0c_r^0\varrho^\alpha$, while \pref{lemest} guarantees that \eqref{thmsol1} holds. Hence, \tref{thmsol} applies and yields \eqref{box1} with $K_\ast:=K_0c_r^0$. In particular, for $N_0\in\N$ and $\rho_0>0$ from \tref{thmsol} there is a $N_1\geq N_0$ so that $\sup_{t\in\Z}\norm{\phi_t^n-\phi_t^0}_\alpha<\tfrac{\rho_0}{2}$ for all $n\geq N_1$. 

	(a) Let $\rho_1>0$ be so small that the sequence $(\phi_t)_{\tau\leq t}$ from \tref{thmbundle}(a) satisfies $\norm{\phi_t}_\alpha<\tfrac{\rho_0}{2}$ for all $\tau\leq t$; such a $\rho_1$ exists since the sequence is contained in the stable fiber bundle of $\phi^0$. Furthermore, for each $\vartheta\in[0,1]$ we obtain
	$$
		\abs{\phi_t^n(y)+\vartheta\phi_t(y)-\phi_t^0(y)}
		\leq
		\norm{\phi_t^n+\vartheta\phi_t-\phi_t^0}_0
		<
		\rho_0\fall y\in\Omega.
	$$
	Now set $\bar r:=r+\rho_0$. Combining the triangle inequality, \cref{corlip} and \pref{lemest} yields that there exist $L_{\bar r}\geq 0$ such that
	\begin{eqnarray*}
		&&
		\norm{[D\sF_t^n(\phi_t^n+\vartheta\phi_t)-D\sF_t^0(\phi_t^0+\vartheta\phi_t)]\phi_t}_\alpha\\
		& \stackrel{\eqref{corlip2}}{\leq} &
		\norm{[D\sF_t^n(\phi_t^n+\vartheta\phi_t)-D\sF_t^0(\phi_t^n+\vartheta\phi_t)]\phi_t}_\alpha+
		L_{\bar r}\norm{\phi_t^n-\phi_t^0}_\alpha\norm{\phi_t}_\alpha\\
		& \stackrel{\eqref{lemest3}}{\leq} &
		\frac{c_0c_{\bar r}^1}{n^\alpha}\norm{\phi_t}_\alpha+
		L_{\bar r}\norm{\phi_t^n-\phi_t^0}_\alpha\norm{\phi_t}_\alpha
		\stackrel{\eqref{thmsol4}}{\leq}
		\frac{c_0c_{\bar r}^1}{n^\alpha}\norm{\phi_t}_\alpha+
		\frac{L_{\bar r}K_\ast c_0}{n^\alpha}\norm{\phi_t}_\alpha
	\end{eqnarray*}
	and with $K_+:=c_0c_{\bar r}^1+L_{\bar r}K_\ast c_0$ we obtain
	$$
		\norm{[D\sF_t^n(\phi_t^n+\vartheta\phi_t)-D\sF_t^0(\phi_t^0+\vartheta\phi_t)]\phi_t}_\alpha
		\leq
		\frac{K_+}{n^\alpha}\sup_{\tau\leq s}\norm{\phi_s}_\alpha
		\fall\vartheta\in[0,1],
	$$
	$n\geq N_1$ and $\tau\leq t$. Hence, the claimed estimate follows from \eqref{thmbundle3}. 

	(b) As in (a), applying \eqref{thmbundle6} rather than \eqref{thmbundle3} leads to the assertion. 
	\qed
\end{proof}
\subsection{Differentiable kernel functions}
Convergence rates improving the consistency order $\alpha\in(0,1]$ obtained in \cref{cornyst} can be expected for integrands in \eqref{ide0} being differentiable in $y\in\Omega$. Here we follow the convention to consider a function on a not necessarily open set $\Omega\subset\R^\kappa$ as differentiable, if it allows a differentiable extension to an open superset of $\Omega$. 

Given $p$-times continuously differentiable functions $u:\Omega\to\R^d$ assume that \eqref{quad} allows a \emph{quadrature} or \emph{cubature error} of the form (see \cite{davis:rabinowitz:84})
\begin{equation}
	\abs{\int_\Omega u(y)\d y-Q^n u}
	\leq
	\frac{c_p}{n^p}
	\sup_{x\in\Omega}\abs{D^{p}u(x)}
	\fall n\in\N
	\label{quadest}
\end{equation}
with constants $c_p\geq 0$.

A smooth framework allows the following improvement of \pref{lemnconv}: 
\begin{proposition}[higher order convergence of $\sF_t^n$]\label{prophc}
	Let $t\in\Z$, $p\in\N$ and $\Omega$ be convex. Suppose the kernel function $f_t:\Omega^2\tm Z_t\to\R^d$ fulfills:
	\begin{itemize}
		\item[(iv)] The partial derivative $D_1f_t:\Omega^2\tm Z_t\to L(\R^\kappa,\R^d)$ exists,

		\item[(v)] both $f_t,D_1f_t$ are of class $C_{(2,3)}^{p}$.
	\end{itemize}
	If \eqref{quad} satisfies \eqref{quadest}, then for every $r>0$ there exists a $\bar c_r^0\geq 0$ such that
	\begin{equation}
		\norm{\sF_t^n(u)-\sF_t^0(u)}_\alpha
		\leq
		\frac{c_p\bar c_r^0}{n^p}\fall n\in\N
		\label{prophest}
	\end{equation}
	and $p$-times continuously differentiable functions $u\in U_t$. 
\end{proposition}
\begin{proof}
	Let $t\in\Z$ and with $u\in U_t$ of class $C^p$ it is convenient to define
	\begin{align*}
		F_t^{(1)}:\Omega^2&\to L(\R^\kappa,\R^d),&
		F_t^{(1)}(x,y)&:=D_1f_t(x,y,u(y)).
	\end{align*}
	The estimate \eqref{prophest} for the $\norm{\cdot}_0$-norm is an immediate consequence of the error estimate \eqref{quadest} and the higher-order chain rule. Let $x,\bar x\in\Omega$ and the Mean Value Theorem \cite[p.~341, Thm.~4.2]{lang:93} gives
	\begin{eqnarray*}
		&&
		[\sF_t^n(u)-\sF_t^0(u)](x)-[\sF_t^n(u)-\sF_t^0(u)](\bar x)\\
		& = &
		\sum_{\eta\in\Omega_n}w_\eta f_t(x,\eta,u(\eta))-\int_\Omega f_t(x,y,u(y))\d y\\
		&&
		-\sum_{\eta\in\Omega_n}w_\eta f_t(\bar x,\eta,u(\eta))+\int_\Omega f_t(\bar x,y,u(y))\d y\\
		& = &
		\int_\Omega\int_0^1D_1f_t(\bar x+\vartheta(x-\bar x),y,u(y))\d\vartheta\d y\intoo{x-\bar x}\\
		&&
		-
		\sum_{\eta\in\Omega_n}w_\eta\int_0^1D_1f_t(\bar x+\vartheta(x-\bar x),\eta,u(\eta))\d\vartheta\intoo{x-\bar x}, 
	\end{eqnarray*}
	from which Fubini's theorem \cite[p.~162, Thm.~8.4]{lang:93} yields
	\begin{align*}
		&
		[\sF_t^n(u)-\sF_t^0(u)](x)-[\sF_t^n(u)-\sF_t^0(u)](\bar x)\\
		= &
		\int_0^1\Bigl(\int_\Omega F_t^{(1)}(\bar x+\vartheta(x-\bar x),y)\d y
		-\!
		\sum_{\eta\in\Omega_n}w_\eta F_t^{(1)}(\bar x+\vartheta(x-\bar x),\eta)
		\Bigr)\d\vartheta\intoo{x-\bar x}
	\end{align*}
	and passing to the norm implies
	\begin{eqnarray*}
		&&
		\abs{[\sF_t^n(u)-\sF_t^0(u)](x)-[\sF_t^n(u)-\sF_t^0(u)](\bar x)}\\
		& \stackrel{\eqref{quadest}}{\leq} &
		\frac{c_p}{n^p}
		\int_0^1\sup_{y\in\Omega}\abs{D_2^{p}F_t^{(1)}(\bar x+\vartheta(x-\bar x),y)}\d\vartheta
		\abs{x-\bar x}\\
		& \leq &
		\frac{c_p}{n^p}(\diam\Omega)^{1-\alpha}
		\sup_{x,y\in\Omega}\abs{D_2^{p}F_t^{(1)}(x,y)}
		\abs{x-\bar x}^\alpha.
	\end{eqnarray*}
	Hence, 
	$
		[\sF_t^n(u)-\sF_t^0(u)]_\alpha
		\leq
		\frac{c_p}{n^p}(\diam\Omega)^{1-\alpha}
		\sup_{x,y\in\Omega}\abs{D_2^{p}F_t^{(1)}(x,y)}
	$
	and if we abbreviate 
	$\bar c_r^0:=\max\set{1,(\diam\Omega)^{1-\alpha}}
		\max_{i=0}^1\sup_{x,y\in\Omega}
		\abs{D_2^{p}F_t^{(i)}(x,y)}$,
	then \eqref{normdef} implies the claimed estimate \eqref{prophest}.
	\qed
\end{proof}

Smooth functions $f_t$ and reference solutions $\phi^0$ allow better convergence rates. Indeed under the assumptions of \pref{propU0}, \ref{propUn} and \ref{lemest}, \ref{prophc}, as well as \cref{corlip} results:
\begin{corollary}[saddle-point structure of \eqref{ide0}, $C^p$-case]
	Let $\Omega\subset\R^\kappa$ be convex. Suppose \href{hyp} holds with $\max\set{2,p}\leq m$ and $\phi^0$ is a weakly hyperbolic $\theta_1$-pe\-ri\-odic solution to \eqref{ide0}. If
	\begin{itemize}
		\item[(iv)] \eqref{quad} is stable, has consistency order $\alpha\in(0,1]$ and satisfies \eqref{quadest}, 

		\item[(v)] the partial derivatives $D_1^kf_t:\Omega^2\tm Z_t\to L_k(\R^\kappa,\R^d)$ exists for $0\leq k\leq p$,

		\item[(vi)] both $f_t$ and $D_1f_t$ are of class $C_{(2,3)}^{p}$, 
	\end{itemize}
	then there exist constants $K_\ast,K_+,K_-\geq 0$ and $N_1\in\N$ such that the following holds for all $n\geq N_1$: The associate weakly hyperbolic and $\theta$-periodic solutions $\phi^n$ to \eqref{iden} satisfy
	\begin{equation}
		\sup_{t\in\Z}\norm{\phi_t^n-\phi_t^0}_\alpha
		\leq
		\frac{K_\ast c_p}{n^p}. 
		\label{box2}
	\end{equation}
	If $\phi^0$ is even hyperbolic, then for each $\tau\in\Z$ one has the estimates
	\begin{enumerate}
		\item $\norm{w_+^n(\tau,v)-w_+^0(\tau,v)}_\alpha\leq
		\tfrac{4K}{1-\beta}\tfrac{K_+}{n^p}(1+\sup_{\tau\leq t}\norm{\phi_t}_\alpha)$ 
for all $p$-times continuously differentiable $v\in B_{\rho_1}(0,R(P_\tau))$, 

		\item $\norm{w_-^n(\tau,v)-w_-^0(\tau,v)}_\alpha\leq
		\tfrac{4K}{1-\beta}\tfrac{K_-}{n^p}(1+\sup_{t\leq\tau}\norm{\phi_t}_\alpha)$ for all $v\in B_{\rho_1}(0,N(P_\tau))$ 
	\end{enumerate}
	with the forward resp.\ backward solution $\phi$ to the IDE \eqref{ide0} from \tref{thmbundle}. 
\end{corollary}
\begin{proof}
	Let $t\in\Z$. Above all, as entire solutions to \eqref{ide0} the functions $\phi_t^0$ are of class $C^p$ due to (v) and \cite[p.~355, Thm.~8.1]{lang:93}. By means of \pref{prophc} the estimate \eqref{box2} results as in the above proof of \cref{cornyst}, with \eqref{lemnconv2} replaced by \eqref{prophest}. 

	(a) As in the above proof of \cref{cornyst} one obtains
	\begin{eqnarray}
		&&
		\norm{[D\sF_t^n(\phi_t^n+\vartheta\phi_t)-D\sF_t^0(\phi_t^0+\vartheta\phi_t)]\phi_t}_\alpha
		\label{tempy}\\
		& \stackrel{\eqref{prophest}}{\leq} &
		\norm{[D\sF_t^n(\phi_t^n+\vartheta\phi_t)-D\sF_t^0(\phi_t^n+\vartheta\phi_t)]\phi_t}_\alpha
		+\frac{L_{\bar r}c_p\bar c_{\bar r}^0}{n^p}\norm{\phi_t}_\alpha
		\notag
	\end{eqnarray}
	for all $\vartheta\in[0,1]$ and $n\geq N_1$. \pref{propU0} yields the explicit derivative 
	$$
		D\sF_t^0(\phi_t^n+\vartheta\phi_t)\phi_t
		\stackrel{\eqref{propDU0a}}{=}
		\int_\Omega D_3f_t(\cdot,y,\phi_t^n(y)+\vartheta\phi_t(y))\phi_t(y)\d y
	$$
	and for the integrand on the right-hand side we observe: Thanks to (v) the periodic solution $\phi^n$ consists of $C^p$-functions $\phi_t^n:\Omega\to\R^d$ and also forward solutions to the IDE \eqref{ide0} are of class $C^p$, i.e. $\phi_t$ is a $C^p$-function for all $t>\tau$. For $t=\tau$ we have $\phi_\tau=v+w_+^0(\tau,v)$ and because $w_+^0(\tau,\cdot)$ is of class $C^m$ by \tref{thmbundle}(a) and $p\leq m$, with $v$ also the initial function $\phi_\tau$ is $p$-times continuously differentiable. Due to (vi) this yields that the integrand $D_3f_t(x,\cdot,\phi_t^n(\cdot)+\vartheta\phi_t(\cdot))\phi_t(\cdot):\Omega\to\R^d$ is of class $C^p$ and the estimate \eqref{quadest} applies. Hence, as in the proof of \pref{prophc} one shows that there exists a $\tilde C\geq 0$ so that
	$
		\norm{[D\sF_t^n(\phi_t^n+\vartheta\phi_t)-D\sF_t^0(\phi_t^n+\vartheta\phi_t)]\phi_t}_\alpha
		\leq
		\tfrac{\tilde C}{n^p}
	$
	and whence \eqref{tempy} yields for all $\vartheta\in[0,1]$, $n\geq N_1$ and $\tau\leq t$ that
	$$
		\norm{[D\sF_t^n(\phi_t^n+\vartheta\phi_t)-D\sF_t^0(\phi_t^0+\vartheta\phi_t)]\phi_t}_\alpha
		\leq
		\intoo{\tilde C+L_rC\sup_{\tau\leq s}\norm{\phi_s}_\alpha}\frac{1}{n^p}.
	$$
	Therefore, the estimate (a) follows from \eqref{thmbundle3}. 

	(b) As above in (a), applying \eqref{thmbundle6} rather than \eqref{thmbundle3} leads to the claimed estimate. Note here that $(\phi_t)_{t\leq\tau}$ is a backward solution to \eqref{ide0} and consequently consists of $C^p$-solutions. Whence, also the initial value $\phi_\tau=v+w_-^0(\tau,v)$ is of class $C^p$ and it is not necessary to assume $v$ to be smooth. 
	\qed
\end{proof}
\begin{appendix}
\renewcommand{\theequation}{\Alph{section}.\arabic{equation}}
\section{Qualitative implicit and Lipschitz inverse function theorem}
\label{appA}
Let $\sX,\sY$ be Banach spaces. We formulate an abstract, but tailor-made implicit function theorem, whose parameter set is merely supposed to be a metric space $(\Lambda,d)$: 
\begin{theorem}[qualitative implicit function theorem]\label{thmifta}
	Let $\varOmega\subseteq\sX$ be nonempty open, $x_0\in\varOmega$, $\lambda_0\in\Lambda$, $y_0\in\sY$, $q\in[0,1)$, and suppose $T:\varOmega\tm\Lambda\to\sY$ satisfies
	\begin{itemize}
		\item[(i')] $T(x_0,\lambda_0)=y_0$, 

		\item[(ii')] the partial derivative $D_1T:\varOmega\tm\Lambda\to L(\sX,\sY)$ exists with $D_1T(x_0,\lambda_0)\in GL(\sX,\sY)$, 

		\item[(iii')] there exist functions $\Gamma_0:\R_+\to\R_+$ and $\Gamma:\R_+^2\to\R_+$ which satisfy $\lim_{\varrho\searrow 0}\Gamma_0(\varrho)=0$, $\lim_{\varrho_1,\varrho_2\searrow 0}\Gamma(\varrho_1,\varrho_2)=0$, such that for all $x\in\varOmega$, $\lambda\in\Lambda$ it holds
		\begin{align}
			\norm{T(x_0,\lambda)-T(x_0,\lambda_0)}
			&\leq
			\Gamma_0(d(\lambda,\lambda_0)),
			\label{thmift1}\\
			\norm{D_1T(x,\lambda)-D_1T(x_0,\lambda_0)}
			&\leq
			\Gamma(\norm{x-x_0},d(\lambda,\lambda_0)).
			\label{thmift2}
		\end{align}
	\end{itemize}
	If $K:=\norm{D_1T(x_0,\lambda_0)^{-1}}$ and $\rho_0,\delta>0$ are chosen so small that
	\begin{align}
		\Gamma_0(\delta)&\leq\tfrac{1-q}{K}\rho_0,&
		\Gamma(\rho_0,\delta)&\leq\tfrac{q}{K},
		\label{thmift3}
	\end{align}
	then there exists a function $\phi:B_\delta(\lambda_0,\Lambda)\to\bar B_{\rho_0}(x_0,\sX)$ satisfying 
	\begin{enumerate}
		\item $\phi(\lambda_0)=x_0$, 

		\item $T(x,\lambda)=y_0$ in $\bar B_{\rho_0}(x_0,\sX)\tm B_\delta(\lambda_0,\Lambda)$ if and only if $x=\phi(\lambda)$, 

		\item $\norm{\phi(\lambda)-x_0}\leq\tfrac{K}{1-q}\Gamma_0(d(\lambda,\lambda_0))$ for all $\lambda\in B_\delta(\lambda_0,\Lambda)$. 
	\end{enumerate}
\end{theorem}
\begin{proof}
	The proof is similar to the one of \cite[Thm.~A.1]{poetzsche:18a}. 
	\qed
\end{proof}

\begin{theorem}[Lipschitz inverse function theorem]\label{thmlip}
	Let $x_0\in\sX$ and $\rho>0$ be given. If a mapping $T:\bar B_\rho(x_0,\sX)\to\sY$ is of the form $T=A+G$ with 
	\begin{enumerate}
		\item[(i)] $A\in GL(\sX,\sY)$, 

		\item[(ii)] $G:\bar B_\rho(x_0,\sX)\to\sY$ is Lipschitz with Lipschitz constant $l<\norm{A^{-1}}^{-1}$, 
	\end{enumerate}
	then the following holds with $\sigma\in\bigl(l,\norm{A^{-1}}^{-1}\bigr]$: 
	\begin{enumerate}
		\item For all $x,\bar x\in\bar B_\rho(x_0,\sX)$ one has
		\begin{equation}
			(\sigma-l)\norm{x-\bar x}
			\leq
			\norm{T(x)-T(\bar x)}
			\leq
			(\norm{A}+l)\norm{x-\bar x},
			\label{thmlip1}
		\end{equation}

		\item for all $y\in\bar B_{(\sigma-l)\rho}(T(x_0),\sY)$ the equation $T(x)=y$ has a unique solution $x^\ast(y)\in\bar B_\rho(x_0,\sX)$, 

		\item with $G|_{B_\rho(x_0,\sX)}$ also the function $x^\ast:B_{(\sigma-l)\rho}(T(x_0),\sY)\to\sX$ is of class $C^m$, $m\in\N_0$.
	\end{enumerate}
\end{theorem}
\begin{proof}
	See \cite[p.~224, (C.11)]{irwin:80}, with the smoothness assertion resulting from the Uniform Contraction Principle \cite[p.~25, Thm.~2.2]{chow:hale:96}.
	\qed
\end{proof}
\end{appendix}
%
%
%
\providecommand{\bysame}{\leavevmode\hbox to3em{\hrulefill}\thinspace}
\providecommand{\MR}{\relax\ifhmode\unskip\space\fi MR }
\providecommand{\MRhref}[2]{%
\href{http://www.ams.org/mathscinet-getitem?mr=#1}{#2}
}
\providecommand{\href}[2]{#2}

\end{document}